\documentclass [11pt,twoside,a4paper]{article}
\usepackage{amsfonts}
\usepackage{amsmath}
\usepackage{amstext}
\usepackage{amssymb}
\usepackage{mathrsfs}
\usepackage{amscd}
\usepackage{xypic}
\usepackage{epsf}              
\usepackage{graphicx}          
\usepackage{fancybox}          
\usepackage{color}             
\usepackage{fancyhdr}
\usepackage[hang,footnotesize]{caption2}  
\usepackage{mathbbold,mathabx}          

\newfont{\aaa}{cmb10 at 19pt}
\newfont{\bbb}{cmb10 at 11pt}
\newtheorem{thm}{Theorem}[section]
\newtheorem{cor}[thm]{Corollary}
\newtheorem{lem}[thm]{Lemma}
\newtheorem{prop}[thm]{Proposition}

\newtheorem{rem}[thm]{Remark}
\newtheorem{exm}[thm]{Example}
\makeatletter
\def\@evenhead{
\vbox{\hbox to \textwidth {}{\hspace{0mm}{\footnotesize
\thepage}}{\hspace{10cm} {\footnotesize {Mu-Fa CHEN et al.}}} \protect\vspace{1truemm}\relax \hrule depth0pt
height0.15truemm width\textwidth}}
\def\@evenfoot{}
\def\@oddhead{\vbox{\hbox to \textwidth
{{\hspace{0cm}{\footnotesize Mixed eigenvalues of discrete $p\,$-Laplacian}\hfill{\footnotesize
\thepage}}\hspace{0mm}}{} \protect\vspace{1truemm}\relax\hrule
depth0pt height0.15truemm width\textwidth}}
\def\@oddfoot{}
\makeatother

\def\scr{\mathscr}
\def\supp{\text{\scriptsize\rm supp\,}}

\newcommand{\rf}[2]{[\ref{#1}; #2]}

\begin{document}

\thispagestyle{empty} \thispagestyle{fancy} {
\fancyhead[lO,RE]{\footnotesize Front. Math. China  2014, 9(6): 1261--1292\\
DOI 10.1007/s11464-013-0186-5\\[3mm]
\includegraphics[0,-50][0,0]{11.bmp}}
\fancyhead[RO,LE]{\scriptsize \bf 
} \fancyfoot[CE,CO]{}}
\renewcommand{\headrulewidth}{0pt}


\setcounter{page}{1}
\qquad\\[8mm]

\noindent{\aaa{Mixed  Eigenvalues of Discrete {\LARGE $\pmb p$\,}-Laplacian }}\\[1mm]

\noindent{\bbb Mu-Fa CHEN$^{1}$,\quad Ling-Di WANG$^{1,2}$,\quad Yu-Hui ZHANG$^{1}$}\\[-1mm]

\noindent\footnotesize{1 School of Mathematical Sciences, Beijing Normal University, Beijing 100875, China\\
2 School of Mathematics and Information Sciences, Henan University, Kaifeng, Henan 475004, China}\\[6mm]

\vskip-2mm \noindent{\footnotesize$\copyright$ Higher Education
Press and Springer-Verlag Berlin Heidelberg 2013} \vskip 4mm

\normalsize\noindent{\bbb Abstract}\quad This paper deals with the principal eigenvalue of discrete $p$-Laplacian on the set
of nonnegative integers. Alternatively, it is studying the optimal constant of a
class of weighted Hardy inequalities. The main goal is the quantitative estimates
of the eigenvalue. The paper begins with the case having reflecting boundary at
origin and absorbing boundary at infinity. Several variational formulas are
presented in different formulation: the difference form, the single summation
form and the double summation form. As their applications, some  explicit lower
and upper  estimates, a criterion for positivity (which was known years ago),
as well as an approximating procedure for the eigenvalue are obtained. Similarly,
the dual case having absorbing boundary at origin and reflecting boundary at
infinity is also studied. Two examples are presented at the end of Section 2
to illustrate the value of the investigation.\vspace{0.3cm}

\footnotetext{Received August 13, 2013; accepted March 3, 2014\\
\hspace*{5.8mm}Corresponding author: Ling-Di WANG, E-mail:wanglingdi@mail.bnu.edu.cn
}

\noindent{\bbb Keywords}\quad Discrete $p$-Laplacian, mixed eigenvalue, variational formula, explicit estimate, positivity criterion, approximating procedure
\smallskip
\\
{\bbb MSC}\quad 60J60, 34L15\\[0.4cm]

\setcounter{section}{1}
\setcounter{thm}{0}

\noindent{\bbb 1\quad Introduction}\\[0.1cm]
\noindent
In the past years, we have been interested in various aspects of stability speed,
such as exponentially ergodic rate, exponential decay rate, algebraic convergence speed,
exponential convergence speeds.  The convergence speeds are often described
by principal  eigenvalues or the optimal constants in different types of inequalities.
Having a great effort on the $L^2$-case (refer to [\ref{Chen4}--\ref{Chen1}] and
reference therein), we now come to a more general setup, studying the nonlinear
$p$-Laplacian, especially on the discrete space
$E:=\{0, 1,\cdots,N\}\, (N\leqslant\infty)$ in this paper. This is a typical topic
in harmonic analysis (cf. \cite{KMP}). The method adopted in this paper is analytic
rather than probabilistic. Let us presume that $N<\infty$ for a moment. Following the
classification given in \cite{Chen1,CWZ}, where $p=2$ was treated,  we have four
types of boundary conditions: DD, DN, ND, and NN, according to Dirichlet (code `D') or Neumann (code `N') boundary at each of the endpoints. For instance, the Neumann condition at the right endpoint means that  $f_{N+1}=f_N$. For Dirichlet condition, it means that $f_{N+1}=0$.
 In the continuous context, the DN-case was partially  studied  in \cite{JM}
 by Jin and Mao. For the NN-  and DD-cases, one may refer to \cite{Chen5}.
 Based on  \cite{Chen1, CWZ}, here, we study the mixed eigenvalues
 (i.e., the ND-  and DN-cases) of discrete $p$-Laplacian. Certainly,  the above
 classification of the boundaries for $p$-Laplacian remains meaningful even if
 $N=\infty$.

 The paper is organized as follows.   In Section 2, we  study the ND-case. First, we
 introduce three groups of variational formulas for the eigenvalue.  As a consequence,
 we obtain the basic estimates (i.e., the ratio of the upper and the lower bounds is
 a constant) of the eigenvalue. Furthermore, an approximating procedure and some
 improved estimates are presented. Except the basic estimates, when $p\neq2$,
 the other results seem to be new.
 To illustrate the power of our main results,
 two examples are included at the end of Section 2.
 Usually, the nonlinear case (here it means $p\neq2$)
 is much harder than the linear one ($p=2$). We are lucky in the present situation
 since most of ideas developed in \cite{Chen1} are still suitable in the present
 general setup. This saves us a lot of spaces. Thus, we do not need to publish
 all details, but emphasize some key points and the difference to \cite{Chen1}.
 The sketched proofs are presented in Section 3.
 In Section 4, the corresponding results for the DN-case are presented.
\noindent\\[4mm]

\setcounter{section}{2}
\setcounter{thm}{0}
\noindent{\bbb 2\quad ND-case}\\[0.1cm]
\noindent
Throughout the paper, denoted by $\scr{C}_K$ the set of functions  having compact support. In this section, let $E=\{i: 0\leqslant i<N+1\}$ $(N\leqslant \infty)$. The discrete $p$-Laplacian is defined as follows:
$$\Omega_pf(k)=\nu_k|f_k-f_{k+1}|^{p-2}\big(f_{k+1}-f_k\big)-\nu_{k-1}|f_{k-1}-f_k|^{p-2}\big(f_k-f_{k-1}\big), \quad p> 1,$$
where $\{\nu_k: k\in E\}$ is a positive sequence with boundary condition $\nu_{-1}=0$ (and $f_{-1}=f_0$). Alternatively, we may rewrite $\Omega_p$ as
$$\Omega_pf(k)=\nu_k|f_k-f_{k+1}|^{p-1}{\rm sgn}\big(f_{k+1}-f_k\big)-\nu_{k-1}|f_{k-1}-f_k|^{p-1}{\rm sgn}\big(f_k-f_{k-1}\big), $$
especially when $p\in (1, 2)$.
 Then we have the following discrete version of the $p$-Laplacian eigenvalue problem with ND-boundary conditions:
\begin{eqnarray}
\text{`Eigenequation':}& \Omega_pg(k)=-\lambda\mu_k|g_k|^{p-2}g_k,\quad k\in E;\label{fr2}\\
\text{ND-boundary  conditions:}& 0\neq g_0= g_{-1} \text{ and } g_{N+1}=0 \text{ if } N<\infty.\label{B1}
\end{eqnarray}
If $(\lambda, g)$ is a solution to the eigenvalue problem, then $\lambda$ is called an ``eigenvalue'' and $g$ is its eigenfunction. Especially, when $p=2$, the first (or principal) eigenvalue corresponds to the exponential decay rate for birth-death process on half line, where $\{\mu_k\}$ is just the invariant measure of the birth-death process and $\{\nu_k\}$ is a quantity related to the recurrence criterion of the process (\rf{Chen1}{Sections 2 and  3}).

Define
$$D_p(f)=\sum_{k\in E}\nu_k|f_k-f_{k+1}|^p,\quad p\geqslant1,\qquad f\in \scr{C}_K,$$
and the ordinary inner product
$$(f, g)=\sum_{k\in E}f_kg_k.$$
Then we have
$$D_p(f)=\big(-\Omega_pf, f\big).$$
Actually,
$$\aligned
\big(-\Omega_pf, f\big)\!=&\!\sum_{k=0}^N\nu_kf_k|f_k-f_{k+1}|^{p-2}\big(f_k-f_{k+1}\big)\\
&+\sum_{k=0}^N\nu_{k-1}f_k|f_k-f_{k-1}|^{p-2}\big(f_k-f_{k-1}\big).
\endaligned$$
Since $\nu_{-1}=0$, one may rewrite the second term as $\sum_{k=1}^N$ and then as $\sum_{k=0}^{N-1}$ by a change of the index. Combining the resulting sum with the first one, we get
$$\aligned
\big(-\Omega_pf, f\big)&=\sum_{k=0}^{N-1}\nu_k|f_k-f_{k+1}|^{p-2}\big(f_k-f_{k+1}\big)^2+\nu_N|f_N|^p\\
&=\sum_{k\in E}\nu_k|f_k-f_{k+1}|^p\quad(\text{since } f_{N+1}\!=\!0).
\endaligned$$

In this section, we are interested in the principal eigenvalue defined by the following classical variational formula:
\begin{equation}\label{NDV}
\lambda_p=\inf\big\{D_p(f):\mu\big(|f|^p\big)=1, f\in \scr{C}_K\big\},
\end{equation}
where  $\mu(f)=\sum_{k\in E}\mu_kf_k$. We mention that the Neumann boundary at  left endpoint is described by $f_0=f_{-1}$ or $\nu_{-1}=0$. The Dirichlet boundary condition at right endpoint is described by  $f_{N+1}=0$ if $N<\infty$. Actually, the condition also holds even if  $N=\infty$ ($f_N:=\lim _{i\to N} f_i$ provided $N=\infty$) as will be proved in Proposition $\ref{Np1}$ below.
Formula \eqref{NDV} can be rewritten as the following weighted Hardy inequality:
$$\mu\big(|f|^p\big)\leqslant A D_p(f),\qquad f\in\scr{C}_K$$
with optimal constant $A\!=\!\lambda_p^{-1}$. This explains  the relationship between the $p$-Laplacian eigenvalue and the Hardy's inequality.
Throughout this paper, we concentrate on $p\in (1, \infty)$ since the degenerated
cases that $p=1$ or $\infty$ are often easier (cf. \rf{B-A}{Lemmas 5.4, 5.6 on Page 49 and 56, respectively}).

\subsection{Main results}
To state our main results, we need some notations.
For $p>1$, let $p^*$ be its conjugate number (i.e., $1/p+1/p^*=1$). Define $\hat{\nu}_j=\nu_j^{1-p^*}$ and  three operators which are parallel to those introduced in  \cite{Chen1}, as follows:
$$\aligned
&I_{i}(f)=\frac{1}{\nu_i(f_{i}-f_{i+1})^{p-1}}\sum_{j=0}^{i}\mu_{j}f_{j}^{p-1}\qquad(\text{single summation form}),\\
&I\!I_{i}(f)=\frac{1}{f_{i}^{p-1}}\bigg[\sum_{j\in{\supp}(f) \cap [i, N]}\hat{\nu}_j\bigg(\sum_{k=0}^{j}\mu_{k}f_{k}^{p-1}\bigg)^{p^*-1}\bigg]^{p-1}\\
&\qquad\qquad\qquad(\text{double summation form}),\\
&R_{i}(w)=\mu_i^{-1}\big[\nu_i(1-w_{i})^{p-1}-\nu_{i-1}\big(w_{i-1}^{-1}-1\big)^{p-1}\big]\qquad(\text{difference form}).\endaligned$$
We make a convention that $w_{-1}>0$ is free and $w_N=0$ if $N<\infty$.
For the lower estimates to be studied below, their domains  are defined, respectively, as follows:
 $$\aligned
&\mathscr{F}_{I}=\{f: f>0 \text{ and } f  \text{ is strictly decreasing}\},\\
&\mathscr{F}_{I\!I}=\{f:\ f>0  \text{ on } E\},\\
&\mathscr{W}=\bigg\{w: w_{i}\in(0,1)\text{ if } \sum_{j\in E}\hat{\nu}_j<\infty\text{ and } w_{i}\in(0,1] \text{ if } \sum_{j\in E}\hat{\nu}_j=\infty\bigg\}.
\endaligned$$
For the upper estimates, some modifications are needed to avoid the non-summable problem:
$$
\aligned
&\widetilde{\mathscr{F}}_{I}=\{f: f \ \text{is strictly decreasing on some}\, [n, m], 0\leqslant n<m<N+1, \\
& \qquad \qquad \;\; f_{\cdot}=f_{\cdot\vee n}
\mathbbold{1}_{.\leqslant m}\},\\
&\widetilde{\mathscr{F}}_{I\!I}=\{f:f_i>0 \text{ up to some } m\in[1,N+1) \text{ and then vanishes}\},\\
&\widetilde{\mathscr{W}}=\Big\{w: \exists\,  m\in[1,N+1)\text{ such that }  w_{i}>0 \text{ up to } m-1,  w_{m}=0, \\&\hskip1.8cm w_{i}<1-{(\nu_{i-1}/\nu_i)^{p-1}}\big(w_{i-1}^{-1}-1\big)\text{ for } i=0,1,\cdots,m\Big\}.
\endaligned
$$
In some extent, these functions are imitated of eigenfunction corresponding to $\lambda_p$.
Each part of Theorem $\ref{Nt1}$ below plays a different role in our study. Operator $I$ is used to deduce the basic estimates (Theorem $\ref{Nt2}$) and operator $I\!I$ is a  tool to produce our  approximating procedure (Theorem $\ref{Nt3}$). In comparing with these two operators, the operator   $R$  is easier in the computation. Noting that for each $f\in\mathscr{F}_I$, the term $\inf_{i\in E}I_i(f)^{-1}$ given in part $(1)$ below  is a lower bound of $\lambda_p$, it indicates that  the formulas on the right-hand side of each term in Theorem $\ref{Nt1}$ are mainly  used for the lower estimates. Similarly, the formulas on the left-hand side are used for the upper estimates.
\begin{thm}\label{Nt1}
For $\lambda_p$ $(p>1)$, we have
\begin{itemize}\setlength{\itemsep}{-0.8ex}
\item[$(1)$] single summation forms:
$$\aligned
\inf_{f\in\widetilde{\mathscr{F}}_{I}}\sup_{i\in E}{I}_{i}(f)^{-1}=\lambda_p=\sup_{f\in\mathscr{F}_{I}}\inf_{i\in E}I_{i}(f)^{-1}.
\endaligned$$
\item[$(2)$] Double summation forms:
$$\aligned\inf_{f\in\widetilde{\mathscr{F}}_{I\!I}}\sup_{i\in {\rm supp}(f)}I\!I_{i}(f)^{-1}=\lambda_p=\sup_{f\in\mathscr{F}_{I\!I}}\inf_{i\in E}I\!I_{i}(f)^{-1}.
\endaligned$$
\item[$(3)$]Difference forms:
$$\aligned
\inf_{w\in\widetilde{\mathscr{W}}}\sup_{i\in E}R_{i}(w)=\lambda_p=\sup_{w\in\mathscr{W}}\inf_{i\in E}R_{i}(w).
\endaligned$$
\end{itemize}
Moreover, the supremum on the right-hand sides of the three above  formulas  can be attained.
\end{thm}

 The next proposition adds some additional sets of test functions for operators $I$ and $I\!I$. For simplicity, in what follows, we use $\downarrow$ (resp. $\downdownarrows$) to denote decreasing (resp. strictly decreasing). In parallel, we also use the notation $\uparrow$ and $\upuparrows$.
\begin{prop}\label{Nprop}For $\lambda_p$ $(p>1)$, we have
$$\aligned
\inf_{f\in\widetilde{\mathscr{F}}_{I}}\sup_{i\in {\rm supp}(f)}I\!I_{i}(f)^{-1}=\lambda_p=\sup_{f\in\mathscr{F}_{I}}\inf_{i\in E}I\!I_{i}(f)^{-1},\endaligned$$
$$\aligned
\lambda_p=\inf_{f\in\widetilde{\mathscr{F}}'_{I\!I}\cup\widetilde{\mathscr{F}}_{I\!I}}\sup_{i\in {\rm supp}(f)}I\!I_{i}(f)^{-1}=\inf_{f\in\widetilde{\mathscr{F}}'_{I}}\sup_{i\in E}I_{i}(f)^{-1},
\endaligned$$
where
$$
\aligned
&\widetilde{\mathscr{F}}'_{I}=\{f: f \downdownarrows,  f \text{ is positive up to some}\ m\in[1,N+1),\; \text{then vanishes}\}\!\subset\! \widetilde{\mathscr{F}}_{I},\\
&\widetilde{\mathscr{F}}'_{I\!I}=\big\{f:f>0 \mbox{ and }\  fI\!I(f)^{p^*-1}\in L^{p}(\mu)\big\}.
\endaligned$$
\end{prop}

Throughout the paper, we write
     ${\tilde\mu}[m, n]=\sum_{j=m}^n {\tilde\mu}_j$ for a measure ${\tilde\mu}$
     and define $k(p)=p{p^*}^{p-1}$ (Figure 1).
     Next, define
     $$\sigma_p=\sup_{n\in E}\big(\mu[0, n]{\hat\nu}[n, N]^{p-1}\big).$$
\begin{figure}[h]
\begin{center}{\includegraphics[width=11.0cm,height=7.5cm]{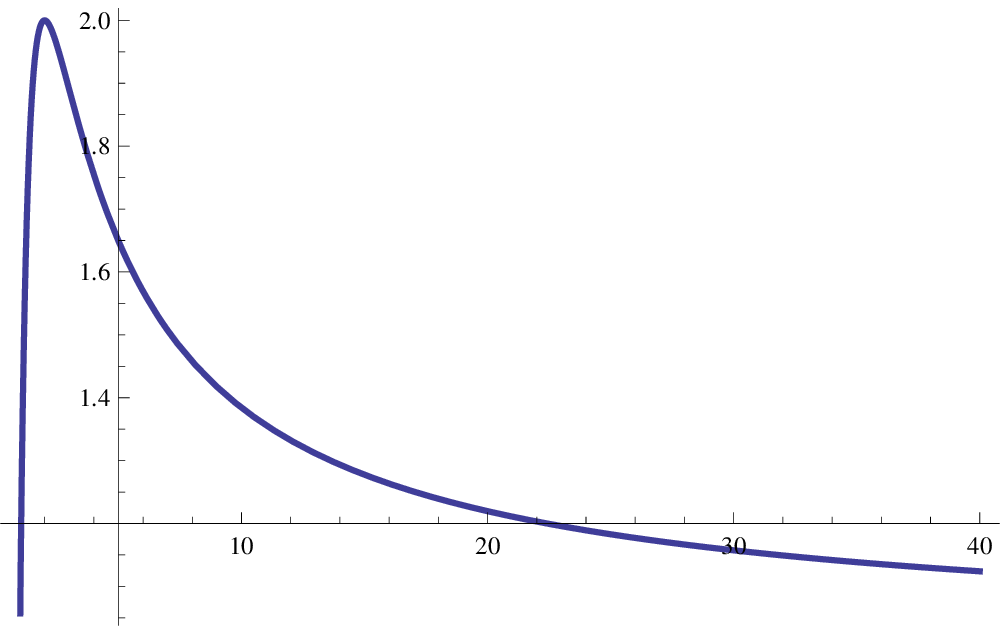}}\end{center}
{\bf Figure 1}\quad Function
$p\to k(p)^{1/p}$ is unimodal with maximum $2$ at $p=2$.
\end{figure}

Applying $f=\hat{\nu}[\cdot, N]^{r(p-1)}$  ($r=1/2$ or $1$) to Theorem \ref{Nt1}\,(1), we obtain the basic estimates given in Theorem $\ref{Nt2}$ below.  This result was known in 1990's (cf. \rf{KMP}{Page 58, Theorem7}). See also \cite{Mao}.
\begin{thm}\label{Nt2}$(\text{\rm Basic estimates})$ For $p>1$, we have $\lambda_p>0$ if and only if $\sigma_p<\infty$.
More precisely,
$$(k(p)\sigma_p)^{-1}\leqslant\lambda_p\leqslant\sigma_p^{-1}.$$
In particular, when $N=\infty$, we have
$$\lambda_p=0\text{ if }\hat{\nu}[1,\infty)=
\infty\text{ and }\lambda_p>0\text{ if }\sum_{k=0}^{\infty}\hat{\nu}_k\mu[0,k]^{p^*-1}<\infty.$$
\end{thm}
 As an application of variational formulas in Theorem \ref{Nt1}\,(2), we obtain an  approximating procedure in the next theorem. This approach can improve the above basic estimates step by step. Noticing that  $\lambda_p$ is trivial once $\sigma_p=\infty$ by Theorem \ref{Nt2}, we may assume that $\sigma_p<\infty$ in the study on the  approximating procedure.
\begin{thm}\label{Nt3}$(\text{\rm Approximating procedure})$ Assume that $\sigma_p<\infty$. Let $p>1$.
\begin{itemize}\setlength{\itemsep}{-0.8ex}
\item[$(1)$] When $\hat{\nu}[0, N]<\infty$, define
    $$f_1=\hat{\nu}[\cdot, N]^{1/p^*},\;\; f_n=f_{n-1}I\!I(f_{n-1})^{p^*-1}\;(n\geqslant 2),$$
    and $\delta_n=\sup_{i\in E}I\!I_i(f_n)$. Otherwise, define $\delta_n=\infty$. Then $\delta_n$ is decreasing in $n$ $($denote its limit by $\delta_{\infty}$$)$ and
    $$\lambda_p\geqslant \delta_{\infty}^{-1}\geqslant\cdots\geqslant\delta_1^{-1}\geqslant\big(k(p)\sigma_p\big)^{-1}.$$
\item[$(2)$] For fixed $\ell, m\in E$, $\ell<m$,
 define
  $$f_1^{(\ell,m)}=\hat{\nu}[\cdot\vee \ell,m]\,\mathbbold{1}_{\leqslant m},\qquad f_n^{(\ell,m)}=f_{n-1}^{(\ell,m)}\big(I\!I\big(f_{n-1}^{(\ell,m)}\big)\big)^{p^*-1}\mathbbold{1}_{\leqslant m},\quad n\geqslant2,$$
where $\mathbbold{1}_{\leqslant m}$ is the indicator of the set $\{0,1,\cdots,m\}$ and then define
$$\delta_n'=\sup_{\ell,m: \ell<m}\min_{i\leqslant m}I\!I_i\big(f_n^{(\ell,m)}\big).$$
Then $\delta_n'$ is increasing in $n$ $($denote its limit by $\delta_{\infty}'$$)$ and $$\sigma_p^{-1}\geqslant\delta_1'^{-1}\geqslant\cdots\geqslant\delta_{\infty}'^{-1}\geqslant\lambda_p.$$
Next, define
$$\bar{\delta}_n=\sup_{\ell,m:\,\ell<m}\frac{\mu\big(\big|f_n^{(\ell,m)}\big|^p\big)}{D_p\big(f_n^{(\ell,m)}\big)},\qquad n\geqslant1.$$ Then $\bar{\delta}_n^{-1}\geqslant\lambda_p$ and $\bar{\delta}_{n+1}\geqslant\delta_n'$  for $n\geqslant1$.
\end{itemize}
\end{thm}

The next result is a consequence of Theorem \ref{Nt3}.
\begin{cor}\label{Nc1}$($\text{\rm Improved estimates}$)$  For $p>1$,
we have
$$\sigma_p^{-1}\geqslant\delta_1'^{-1}\geqslant\lambda_p\geqslant\delta_1^{-1}\geqslant\big(k(p)\sigma_p\big)^{-1},$$
where
$$\aligned&\delta_1=\sup_{i\in E}\bigg[\frac{1}{\hat{\nu}[i,N]^{1/p^*}}
\sum_{j=i}^N\hat{\nu}_j\bigg(\sum_{k=0}^j\mu_k\hat{\nu}[k, N]^{(p-1)/p^*}\bigg)^{p^*-1}\bigg]^{p-1},\\
&\delta_1'=\sup_{\ell\in E}\frac{1}{\hat{\nu}[\ell, N]^{p-1}}\bigg[\sum_{j=\ell}^N\hat{\nu}_j\bigg(\sum_{k=0}^{j}\mu_{k}\hat{\nu}[k\vee \ell, N]^{p-1}\bigg)^{p^*-1}\bigg]^{p-1}.
\endaligned$$
Moreover,
$$\bar{\delta}_{1}=\sup_{m\in E}\frac{1}{\hat{\nu}[m, N]}\sum_{j=0}^{N}\mu_{j}\hat{\nu}[j \vee m, N]^p\in[\sigma_p, p\sigma_p],$$
and
$\bar\delta_1\leqslant\delta_1'$ for $1<p\leqslant2$, $\bar\delta_1\geqslant\delta_1'$ for $p\geqslant2$.
\end{cor}

An remarkable point of Corollary \ref{Nc1} is its last assertion which is comparable
with the known result that $\bar\delta_1=\delta_1'$ when $p=2$
(cf. \rf{Chen1}{Theorem 3.2}). This indicates that some additional work is
necessary for general $p$ than the specific one $p=2$.
\noindent\\[4mm]

\subsection{Examples}
In the worst case that $p=2$ (cf. Figure 1), the ratio $k(p)^{1/p}$ of the upper and
lower estimates is no more than 2 which can be improved (no more than $\sqrt{2}$\,) by the improved estimates as shown by a large number of examples (cf. \cite{Chen1}). The same conclusion should also be true for general $p$ as shown by two examples below.
Actually, the effectiveness of the improved bounds $\delta_1$ and ${\bar\delta}_1$
shown by the examples is quite unexpected.

\begin{exm}\label{EXTH1}{\rm Assume that $E=\{0, 1,\cdots, N\}$, $a>0$, and $r>1$. Let $\mu_k=r^{k}$, $\nu_k=ar^{k+1}$  for $k\in E$. Then
$$\aligned
\sigma_p&=\frac{r}{a(r-1)\big(r^{p^*-1}-1\big)^{p-1}},\\
\delta_1&=\frac{1}{ar\big(r^{1/p}-1\big)}\sup_{i\in E}\bigg\{\sum_{j=i}^N
  \Big(r^{i-j+[(j-i+1)/p]}-r^{i-j-(i/p)}\Big)^{p^*-1}\bigg\}^{p-1},\\
\bar{\delta}_1&=\frac{r^{p^*}-1}{a\big(r^{p^*-1}-1\big)^p(r-1)},\\
\delta_1'&=\frac{1}{ar}\sup_{\ell\in E}\bigg\{\sum_{j=\ell}^N
\bigg(\frac{r^{\ell+1}-1}{r^j(r-1)}+(j-\ell)r^{\ell-j}\bigg)^{p^*-1}\bigg\}^{p-1}.
\endaligned$$
The improved estimates given in Corollary {\rm\ref{Nc1}} are shown in Figure $2$ below.
}\end{exm}

The ratio between $\delta_1$ and $\delta_1'$ (or ${\bar\delta}_1$) is obvious smaller than the
basic estimates $k(p)$ obtained in Theorem $\ref{Nt2}$. When $p=2$,
${\bar\delta}_1=\delta_1'$ which is known as just mentioned.
 \begin{center}{\includegraphics[width=11.0cm,height=7.5cm]{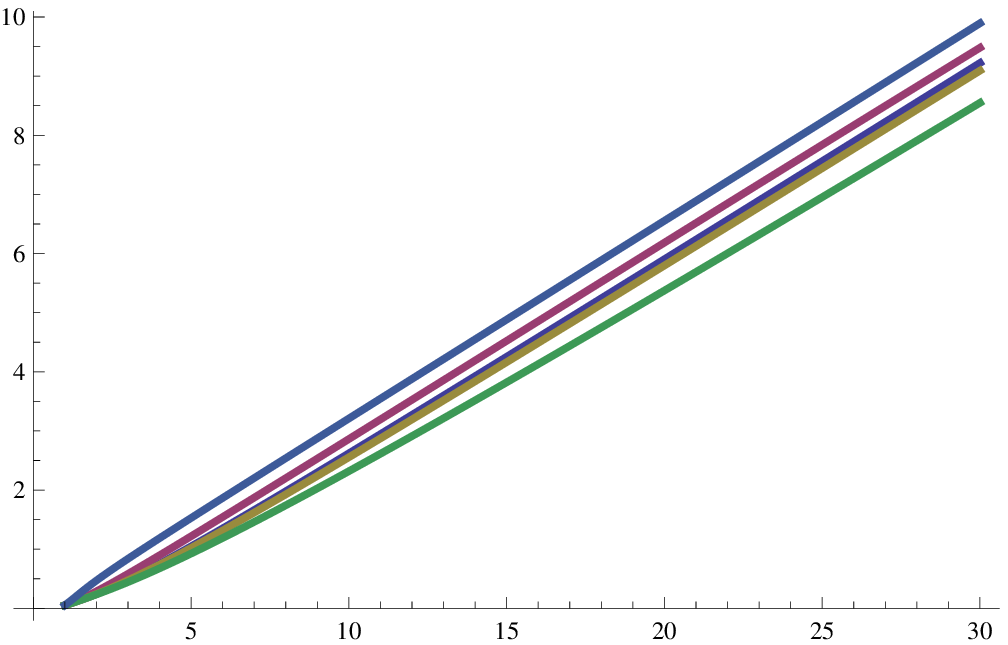}
}\end{center}
{\bf Figure 2}\quad Let $N=80$, $a=1$, $r=20$ and let $p$ vary over $(1.001, 30.001)$
avoiding the singularity at $p=1$.
Viewing from the right-hand side,
the curves from top to bottom are
$(k(p) \sigma_p)^{1/p}$, ${\delta_1}^{1/p}$, ${\bar\delta}_1^{1/p}$, ${\delta_1'}^{1/p}$ and $\sigma_p^{1/p}$, respectively.

Note that the lower bounds ${\bar\delta}_1^{1/p}$ and ${\delta_1'}^{1/p}$ of $\lambda_p^{-1/p}$ are nearly overlapped.

\begin{exm}\label{EXTH4}{\rm Assume that $E=\{0, 1,\cdots,N\}$, $N<\infty$.
Let $\mu_k=1$ and $\nu_k=1$ for $k\in E$. Then
$$\aligned
\sigma_p&=\sup_{n\in E}\big[ (n+1)(N-n+1)^{p-1}\big],\\
\delta_1&=\sup_{i\in E}\bigg[\frac{1}{(N-i+1)^{(p-1)/p}}\sum_{j=i}^N
\bigg(\sum_{k=0}^j\big(N-k+1\big)^{(p-1)/p^*}\bigg)^{p^*-1}\bigg]^{p-1},\\
\bar{\delta}_1&=\sup_{m\in E}\bigg (m(N-m+1)^{p-1}+\frac{1}{N-m+1}\sum_{j=m}^N(N-j+1)^p\bigg),\\
\delta_1'&=\sup_{\ell\in E}\bigg\{\frac{1}{N-\ell+1}\sum_{j=\ell}^N
\bigg[\ell(N-\ell+1)^{p-1}+\sum_{k=\ell}^j\big(N-k+1\big)^{p-1}\bigg]^{p^*-1}\bigg\}^{p-1}.
\endaligned$$
Surprisingly,  the improved estimates $\delta_1$, $\delta_1'$ and $\bar\delta_1$
are nearly overlapped as shown in Figure $3$.
}\end{exm}
\begin{center}{\includegraphics[width=11.0cm,height=7.5cm]{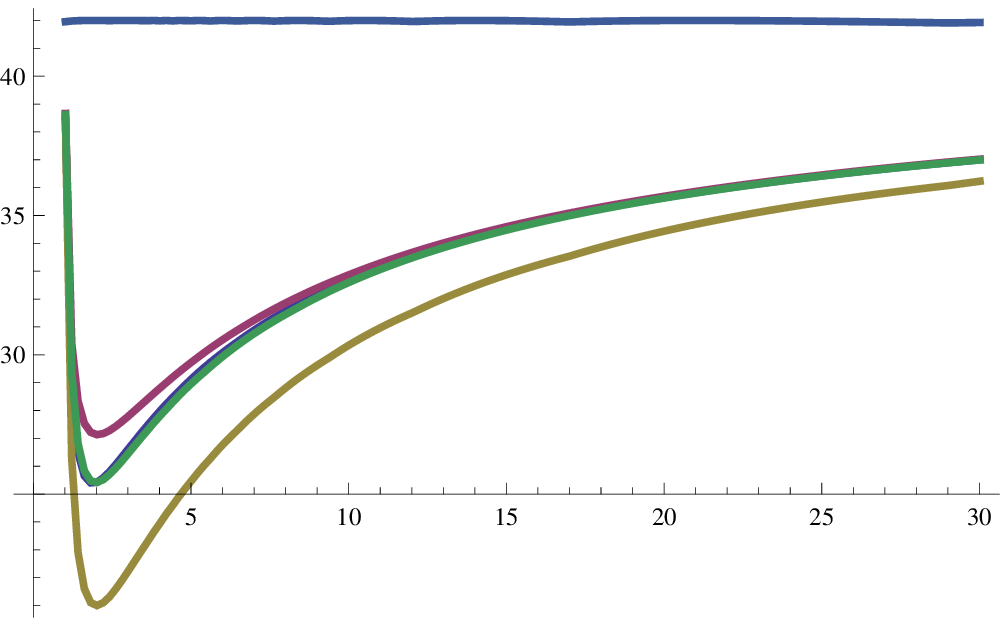}}\end{center}
{\bf Figure 3}\quad Let $N=40$ and let $p$ vary over $(1.0175, 30.0175)$ avoiding
the singularity at $p=1$.
Viewing from the right-hand side,
the curves from top to bottom are again
$(k(p)\sigma_p)^{1/p}$, ${\delta_1}^{1/p}$, ${\bar\delta}_1^{1/p}$, ${\delta_1'}^{1/p}$, and $\sigma_p^{1/p}$, respectively. Note that ${\bar\delta}_1^{1/p}$ and ${\delta_1'}^{1/p}$ (lower bounds of $\lambda_p^{-1/p}$), as well as ${\delta_1}^{1/p}$ (upper bound) are nearly overlapped, except in a small neighborhood of $p=2$.

For this example, the exact $\lambda_p$ is unknown except that
$\lambda_p=\sin^2 \frac{\pi}{2(N+2)}$ when  $p=2.$
\noindent\\[4mm]

\setcounter{section}{3}
\setcounter{thm}{0}

\noindent{\bbb 3\quad {Proofs of the main results in Section 2}}\\[0.1cm]
This section is organized as follows. Some preparations are collected in Subsection 3.1. The preparations may not be needed completely for our proofs here,
but they are useful for the study in a more general setup. The proofs of the main
results are presented in Subsection 3.2.
\noindent\\[4mm]

\noindent{\bbb 3.1\quad Some preparations}\\[0.1cm]
A large part of the results stated in Section 2 depend on the properties of the
 eigenfunction  $g$ of $\lambda_p$. The goal of this subsection is studying these
 properties.

  Define an operator
  $$\widebar\Omega_pf(k)=\Omega_pf(k)-\mu_kd_k|f_k|^{p-2}f_k,\qquad k\in  E,\; p> 1,$$
 where $\{d_k\}_{k\in E}$ is a fixed nonnegative sequence. Then there is  an extended equation of \eqref{fr2}:
  \begin{equation}\label{fr4} \widebar\Omega_pf(k)=-\bar \lambda\mu_k|f_k|^{p-2}f_k,\qquad k\in E,\end{equation}which coincides with equation \eqref{fr2} for $\bar \lambda=\lambda $ if $d_k=0$ for every $k\in E$.
\begin{prop} \label{Np2}
 Define
 $${\widebar D}_p(f)=D_p(f)+\sum_{k\in E}d_k\mu_k|f_k|^p,\qquad f\in \scr{C}_K.$$
 Let
 \begin{equation}\label{fr5}\bar \lambda=\inf
 \big\{\widebar D_p(f):\mu\big(|f|^{p}\big)=1,  f\in\scr{C}_K \text{ and } f_{N+1}=0 \text{ if }N<\infty\big\}.
 \end{equation}
 Then the solution, say $g$, to equation \eqref{fr4} with ND-boundaries  is either positive or negative. In particular, the assertion holds for the eigenfunction of $\lambda_p$.
 \end{prop}
 {\noindent \it Proof }\quad\rm
Since $g_{-1}=g_0$, by making summation from $0$ to $i\in E$ with respect to $k$ on  both sides of \eqref{fr4}, we get \begin{equation}\label{ff1}\nu_i|g_i-g_{i+1}|^{p-2}\big(g_i-g_{i+1}\big)= \sum_{k=0}^{i}\big(\bar \lambda-d_k\big)\mu_k|g_k|^{p-2}g_k,\qquad i\in E.\end{equation}
 If $\bar \lambda=0$,
then the assertion is obvious by \eqref{ff1} and induction. If $\bar \lambda>0$, then $g_0\neq0$ (otherwise $g\equiv0$). Without loss of generality, assume that  $g_{0}=1$ (if not, replace $g$ with $g/g_0$). Suppose that  there exists $k_0$, $1\leqslant k_0< N$  such that  $g_{i}>0$ for $i<k_0$ and $g_{k_0}\leqslant 0$.
Let
$$f_i=g_i\mathbbold{1}_{i< k_0}+\varepsilon\mathbbold{1}_{i=k_0}$$
for some $0<\varepsilon<g_{k_0-1}$. Then $f$ belongs to the setting defining $\bar\lambda$ (cf. \eqref{fr5}).
 Since $g_{k_0}\leqslant0<\varepsilon<g_{k_0-1}$ and $|\varepsilon-g_{k_0-1}|<|g_{k_0}-g_{k_0-1}|$, we have
$$\begin{aligned}
\widebar\Omega_{p}f(k_0\!-\!1)\!&=\widebar\Omega_{p}g(k_0-1)+\nu_{k_0-1}|g_{k_0-1}-g_{k_0}|^{p-2}\big(g_{k_0-1}-g_{k_0}\big)\\
&\hskip1cm-\nu_{k_0-1}|\varepsilon-g_{k_0-1}|^{p-2}\big(g_{k_0-1}-\varepsilon\big)\\
&\geqslant \widebar\Omega_{p}g(k_0-1)-\nu_{k_0-1}|\varepsilon-g_{k_0-1}|^{p-2}\big[\big(g_{k_0-1}-\varepsilon\big)-\big(g_{k_0-1}-g_{k_0}\big)\big]\\
&=\widebar\Omega_{p}g(k_0-1)-\nu_{k_0-1}\big(g_{k_0-1}-\varepsilon\big)^{p-2}\big(g_{k_0}-\varepsilon\big)\\
&>\widebar\Omega_{p}g(k_0-1),
\\
\widebar\Omega_{p}f(k_0)&=-\nu_{k_0}\varepsilon^{p-1}+\nu_{k_0-1}|g_{k_0-1}-\varepsilon|^{p-2}\big(g_{k_0-1}-\varepsilon\big)-\mu_{k_0}d_{k_0}\varepsilon^{p-1}.
\end{aligned}$$
Hence,
$$\begin{aligned}
\widebar D_p(f)&=\big(-\widebar\Omega_{p}f,f\big)\\
&=-\sum_{i=0}^{k_0-2}f_{i}\widebar\Omega_{p}f(i)-f_{k_0-1}\widebar\Omega_{p}f(k_0-1)-\varepsilon\widebar\Omega_{p}f(k_0)\\
&<-\sum_{i=0}^{k_0-2}g_{i}\widebar\Omega_{p}g(i)-g_{k_0-1}\widebar\Omega_{p}g(k_0-1)-\varepsilon\widebar\Omega_{p}f(k_0)\\
&=\bar \lambda\sum_{i=0}^{k_0-1}\mu_{i}|g_{i}|^{p-2}g_{i}^{2}+\varepsilon\big[(\nu_{k_0}+\mu_{k_0}d_{k_0})\varepsilon^{p-1}-\nu_{k_0-1}\big(g_{{k_0}-1}-\varepsilon\big)^{p-1}\big].
\end{aligned}$$
In the second equality, we have used the fact that
\begin{align}\label{ff2}
\nu_i|g_i-g_{i+1}|^{p-2}\big(g_i-g_{i+1}\big)g_i= \sum_{k=0}^{i}(\bar \lambda-d_k)\mu_k|g_k|^p-\sum_{k=0}^{i-1}\nu_k
& |g_k-g_{k+1}|^p,\nonumber\\
& i\in E,
\end{align}
which can be  obtained  from \eqref{fr4},  by a computation similar to that of $\big(-\Omega_pf, f\big)$ given above \eqref{NDV}.
Noticing that \begin{equation*}
\mu\big(|f|^p\big)=\sum_{i=0}^{k_0-1}\mu_{i}|g_{i}|^{p-2}g_{i}^{2}+\mu_{k_0}\varepsilon^{p}\end{equation*}
and
$$\nu_{k_0}+\mu_{k_0}d_{k_0}-\bar \lambda\mu_{k_0}<\nu_{k_0-1}\bigg(\frac{g_{k_0-1}}{\varepsilon}-1\bigg)^{p-1}$$
 for small enough $\varepsilon$, we obtain a contradiction to \eqref{fr5}:
 $$\frac{\widebar D_p(f)}{\mu\big(|f|^p\big)}<\bar\lambda\leqslant \frac{\widebar D_p(f)}{\mu\big(|f|^p\big)}.$$
 This proves the   first assertion and then the second one is obvious.\qquad$\Box$

 Before moving further, we introduce an equation which is somehow more general than eigenequation:
 \begin{equation}\text{Poisson equation: }\qquad\label{fr0}\Omega_pg(k)=-\mu_k|f_k|^{p-2}f_k.\end{equation}
 By putting $f=\lambda_p^{p^*-1}g$, we return to eigenequation.
From \eqref{fr0},  for $i,j \in  E$ with $i<j$  we obtain
\begin{equation}\label{N-D}
\nu_j|g_j-g_{j+1}|^{p-2}\big(g_j-g_{j+1}\big)-\nu_{i-1}|g_{i-1}-g_i|^{p-2}
\big(g_{i-1}-g_i\big)\!=\!\sum_{k=i}^j\mu_k|f_k|^{p-2}f_k.\end{equation}
Moreover, if $g$ is positive and  decreasing, then
\begin{equation}
\label{Nf1}g_{n}-g_{N+1}=\sum_{j=n}^{N}\bigg(\frac{1}{\nu_j}\sum_{k=0}^{j}\mu_k|f_k|^{p-2}f_k\bigg)^{p^*-1},\ \qquad n\in E.
\end{equation}

  Besides Proposition \ref{Np2}, two more propositions are needed. One describes the monotonicity  of the eigenfunction presented in the next proposition, and the other one is about the vanishing property to be presented in Proposition \ref{Np1}.
\begin{prop}\label{Np3}Assume that $(\lambda_p, g)$ is a solution to \eqref{fr2} with ND-boundaries and $\lambda_p>0$. Then the eigenfunction $g$ is strictly monotone. Furthermore,
\begin{equation}\label{Nf3}
\frac{1}{\lambda_p}=\bigg[\frac{1}{g_{n}-g_{N+1}}\sum_{k=n}^{N}\bigg(\frac{1}{\nu_k}\sum_{i=0}^{k}\mu_{i}g_{i}^{p-1}\bigg)^{p^*-1}\bigg]^{p-1},\ \qquad n\in E.
\end{equation}
\end{prop}
{\noindent \it Proof }\quad Without loss of generality, assume that $g_{0}=1$. The first assertion follows by letting $i=0$ and $f=\lambda_p^{p^*-1}g$ in \eqref{N-D}. Moreover, it is clear that $g$ is strictly decreasing. Formula \eqref{Nf3}  then follows from \eqref{Nf1} by letting $f=\lambda_p^{p^*-1}g$.$\qquad\Box$
\medskip

As mentioned above, with $f=\lambda_p^{p^*-1}g$, \eqref{N-D} and \eqref{Nf1} are simple variants of eigenequation \eqref{fr2}. However, for general test function $f$, the left-hand side of the function $g$ defined by \eqref{Nf1} may be far away from the eigenfunction of $\lambda_p$. Nevertheless, we regard the resulting  function (assuming $g_{N+1}=0$) as a mimic of the eigenfunction. This explains where the operator $I\!I$ comes from: it is regarded as an approximation of $\lambda_p^{-1}$ since $I\!I(\lambda_p^{p^*-1} g)\equiv\lambda_p^{-1}$. Next, write $I\!I_i(f)$ as $u_i/v_i$. Then $I_i(f)$ is defined by
$$I_i(f)=\frac{u_i-u_{i+1}}{v_i-v_{i+1}}.$$
In other words, the operator $I$ comes from $I\!I$ in the use of proportional property. The operator $R$ also comes from the eigenequation by setting $w_i=g_{i+1}/g_i$.

\begin{rem}\label{NR1}
{\rm Define
\begin{equation}\label{fr3}\tilde\lambda_p=\inf\big\{D_p(f): \mu\big(|f|^{p}\big)=1\text{ and } f_{N+1}=0\big\}.\end{equation}
\begin{itemize}\setlength{\itemsep}{-0.8ex}
\item[$(1)$] It is easy to check that the assertions in Propositions {\rm\ref{Np2}}
and {\rm\ref{Np3}} also hold for $\tilde\lambda_p$ defined by \eqref{fr3}.
\item[$(2)$] Define
 $$\lambda_p^{(n)}=\inf\big\{D_p(f): \mu\big(|f|^p\big)=1, f=f_{\leqslant n}\big\}.$$ We have $\lambda_p=\tilde\lambda_p$ and $\lambda_p^{(n)}\downarrow \lambda_p$ as $n\uparrow N$. Indeed, it is clear that $$\lambda_p^{(n)}\geqslant\lambda_p\geqslant\tilde\lambda_p.$$ By definition of $\tilde\lambda_p$, for any fixed $\varepsilon>0$, there exists $\bar f\in L^p(\mu)$ such that  $f_{N+1}=0$ and
$$\frac{D_p(\bar f)}{\mu\big(|\bar f|^p\big)}\leqslant\tilde\lambda_p+\varepsilon.$$
Define $f^{(n)}=\bar f\mathbbold{1}_{[0, n)}$. Then $f^{(n)}\in L^p(\mu)$ and
$$D_p\big(f^{(n)}\big)\;\uparrow\; D_p(\bar f),\qquad \mu\big(|f^{(n)}|^p\big)\;\uparrow\; \mu\big(|\bar f|^p\big)\qquad\text{ as } n\to N.$$
Next,   for large enough $n$, we have
$$\lambda_p\leqslant\lambda_p^{(n)}\leqslant\frac{D_p\big(f^{(n)}\big)}
{\mu\big(|f^{(n)}|^p\big)}\leqslant\frac{D_p(\bar f)}{\mu\big(|\bar f|^p\big)}+\varepsilon\leqslant\tilde\lambda_p+2\varepsilon\leqslant \lambda_p+2\varepsilon.$$
By letting $n\to N$ first and then $\varepsilon\downarrow 0$, it follows that $\lambda_p=\tilde\lambda_p$. Actually, we have $\lambda_p^{(n)}\downarrow \lambda_p$.
 \item[$(3)$] The test functions in \eqref{fr3} are described by $f_{N+1}=0$, which can be seen as the imitations of eigenfunction, so the assertion in item $(2)$ above also implies the vanishing property of eigenfunction to some extent.
\end{itemize}
}\end{rem}

 From now on, in this section,  without loss of generality, we assume that {\bf the eigenfunction $g$ corresponding to $\lambda_p$ (or  $\tilde\lambda_p$) is nonnegative and strictly decreasing.}
The monotone property is important to our study. For example, it guarantees the meaning of the operator $I$ defined in Section 2  since the denominator is $f_i-f_{i+1}$ there.
\noindent\\[4mm]

\noindent{\bbb 3.2\quad Sketch proof of the main results in ND-case}\\[0.1cm]
Since a large part of proofs  are  analogies of the case that $p=2$, we need only  to show some keys and some difference between the general $p$ and the specific $p=2$.

\noindent{\it Proof of Theorem $\ref{Nt1}$} \quad We adopt the following circle arguments for the lower estimates:
$$\aligned
\lambda_p&\!\geqslant\!\!\!\sup_{f\in\mathscr{F}_{I\!I}}\inf_{i\in E}I\!I_{i}(f)^{-1}\\
&\!\!=\!\!\sup_{f\in\mathscr{F}_{I}}\inf_{i\in E}I\!I_{i}(f)^{-1}\\
&\!\!=\!\!\sup_{f\in\mathscr{F}_{I}}\inf_{i\in E}I_{i}(f)^{-1}\\
&\!\!\geqslant\!\!\sup_{w\in\mathscr{W}}\inf_{i\in E}R_{i}(w)\\
&\!\!\geqslant\!\!\lambda_p.\endaligned$$

\noindent{\bf Step 1} \quad Prove that  $${\lambda}_p\geqslant\sup_{f\in\mathscr{F}_{I\!I}}\inf_{i\in E}I\!I_{i}(f)^{-1}.$$

Clearly, we have $\lambda_p\geqslant\tilde\lambda_p$ by Remark \ref{NR1} (2).  Let $\{h_k:k\in E\}$ be a positive sequence and let $g$ satisfy $\mu\big(|g|^{p}\big)=1$ and $g_{N+1}=0$. Then
$$\begin{aligned}
1&=\mu\big(|g|^{p}\big)\\
&\leqslant\sum_{i=0}^{N}\mu_{i}\bigg(\sum_{k=i}^{N}\big|g_{k}-g_{k+1}\big|\bigg)^{p}\quad(\text{since } g_{N+1}=0)\\
&=\sum_{i=0}^{N}\mu_{i}\bigg(\sum_{k=i}^{N}{\big|g_{k}-g_{k+1}\big|}\bigg(\frac{\nu_k}{h_{k}}\bigg)^{1/p}\bigg(\frac{h_{k}}{\nu_k}\bigg)^{1/p}\bigg)^{p}\\
&\leqslant\sum_{i=0}^{N}\mu_{i}\bigg[\bigg(\sum_{k=i}^{N}\big|g_{k}-g_{k+1}\big|^{p}\frac{\nu_k}{h_{k}}\bigg)^{1/p}\bigg(\sum_{j=i}^{N}\bigg(\frac{h_j}{\nu_j}\bigg)^{{p^*}/{p}}\bigg)^{1/p^*}\bigg]^{p}\\
&\qquad\qquad\qquad\qquad\qquad\qquad(\text{by H\"{o}lder's inequality})\\&=\sum_{k=0}^{N}\frac{\nu_k}{h_{k}}\big|g_{k}-g_{k+1}\big|^{p}\sum_{i=0}^{k}\mu_{i}\bigg(\sum_{j=i}^{N}\hat{\nu}_jh_{j}^{p^*-1}\bigg)^{p-1}\\
&\qquad\qquad\qquad\qquad(\text{by exchanging the order of the sums})
\\&\leqslant D_p(g)\sup_{k\in E}H_{k},\end{aligned}$$
where
$$H_{k}=\frac{1}{h_{k}}\sum_{i=0}^{k}\mu_{i}\bigg(\sum_{j=i}^{N}\hat{\nu}_jh_{j}^{p^*-1}\bigg)^{p-1}.$$
  For every $f\in\mathscr{F}_{I\!I}$ with $\sup_{i\in E}I\!I_{i}(f)<\infty$, let $$h_{k}=\sum_{j=0}^{k}\mu_{j}f_{j}^{p-1}.$$
  Then
  $$\sup_{k\in E}H_{k}\leqslant\sup_{k\in E}I\!I_{k}(f)$$
by the proportional property.
Hence,
$$D_p(g)\geqslant\inf_{k\in E}I\!I_k(f)^{-1}$$
for every $g$ with $\mu\big(|g|^p\big)=1$,\; $g_{N+1}=0$, and $f\in\mathscr{F}_{I\!I}.$
By making the supremum with respect to $f\in\mathscr{F}_{I\!I}$ first and then the infimum with respect to  $g$ with $g_{N+1}=0$ and $\mu\big(|g|^p\big)=1$, it follows
 $$\lambda_p\geqslant\tilde\lambda_p\geqslant\sup_{f\in\mathscr{F}_{I\!I}}\inf_{i\in E}I\!I_{i}(f)^{-1}.$$

\noindent{\bf Step 2} \quad Prove that
 $$\sup_{f\in\mathscr{F}_{I\!I}}\inf_{i\in E}I\!I_{i}(f)^{-1}=\sup_{f\in\mathscr{F}_{I}}\inf_{i\in E}I\!I_{i}(f)^{-1}=\sup_{f\in\mathscr{F}_{I}}\inf_{i\in E}I_{i}(f)^{-1}.$$
 Using the proportional property, on the one hand, for any fixed $f\in\mathscr{F}_I$, we have
 $$\begin{aligned}\sup_{i\in E}I\!I_{i}(f)\!&=\sup_{i\in E}\bigg[\frac{1}{f_{i}}\sum_{j=i}^{N}\hat{\nu}_j\bigg(\sum_{k=0}^{j}\mu_{k}f_{k}^{p-1}\bigg)^{p^*-1}\bigg]^{p-1}\\
&=\sup_{i\in E}\bigg\{\bigg[\sum_{j=i}^{N}\hat{\nu}_j\bigg(\sum_{k=0}^{j}\mu_{k}f_{k}^{p-1}\bigg)^{p^*-1}\bigg]\bigg/\bigg[\sum_{j=i}^{N}\big(f_{j}\!-\!f_{j+1}\big)\!+\!f_{N+1}\bigg]\bigg\}^{p-1}\\
&\leqslant\sup_{j\in E}\frac{1}{\nu_j(f_{j}-f_{j+1})^{p-1}}\sum_{k=0}^{j}\mu_{k}f_{k}^{p-1}\qquad\big(\text{since }\hat{\nu}_j={\nu_j}^{1-p^*}\big)\\
&=\sup_{j\in E}I_{i}(f).
\end{aligned}$$
Since $\scr{F}_I\subset\scr{F}_{I\!I}$, by making the infimum with respect to $f\in\scr{F}_I$ on  both sides of the inequality above, we have
 $$\sup_{f\in\mathscr{F}_{I\!I}}\inf_{i\in E}I\!I_{i}(f)^{-1}\geqslant\sup_{f\in\mathscr{F}_{I}}\inf_{i\in E}I\!I_{i}(f)^{-1}\geqslant\sup_{f\in\mathscr{F}_{I}}\inf_{i\in E}I_{i}(f)^{-1}.$$
On the other hand,
 for any fixed $f\in\mathscr{F}_{I\!I}$, let $$g=fI\!I(f)^{p^*-1}\in\mathscr{F}_I.$$
 Similar to the proof above, we have
 $$\sup_{i\in E}I_i(g)\leqslant \sup_{i\in E}I\!I_i(f)$$
  Therefore
  $$\sup_{f\in\mathscr{F}_{I}}\inf_{i\in E}I_{i}(f)^{-1}\geqslant\sup_{f\in\mathscr{F}_{I\!I}}\inf_{i\in E}I\!I_{i}(f)^{-1}$$
 and the required assertion holds.

\noindent{\bf Step 3}\quad  Prove that
 $$\sup_{f\in\mathscr{F}_{I\!I}}\inf_{i\in E}I\!I_{i}(f)^{-1}\geqslant\sup_{w\in\mathscr{W}}\inf_{i\in E}R_{i}(w).$$

 First, we  change the form of $R_i(w)$.
Given $w\in\mathscr{W}$, let $u_{i+1}=w_{i}w_{i-1}\cdots w_0$ for $i\geqslant0$, $u_0=1$. Then $u$ is positive, strictly decreasing and $w_i=u_{i+1}/u_i$. By a simple rearrangement, we get
$$\aligned
R_i(w)=\frac{1}{\mu_i}\bigg[\nu_i\bigg(1-\frac{u_{i+1}}{u_i}\bigg)^{p-1}-\nu_{i-1}\bigg(\frac{u_{i-1}}{u_i}-1\bigg)^{p-1}\bigg]=-\frac{1}{\mu_iu_i^{p-1}}\Omega_pu(i).
\endaligned$$

Next, we prove the main assertion. Without loss of generality, assume that $\inf_{i\in E}R_{i}(w)>0$. Let $u$ be the function constructed above and let $f=uR(w)^{p^*-1}>0$. Then $\Omega_pu=-\mu f^{p-1}$. Since $\nu_{-1}=0$, $u$ is decreasing and $u>0$, by \eqref{Nf1}, we have $$u_{i}-u_{N+1}=\sum_{k=i}^{N}\bigg(\frac{1}{\nu_k}\sum_{j=0}^{k}\mu_{j}f_{j}^{p-1}\bigg)^{p^*-1}.$$ Hence,  $$R_i(w)^{1-p^*}=\frac{u_{i}}{f_{i}}\geqslant\frac{1}{f_i}\sum_{k=i}^{N}\bigg(\frac{1}{\nu_k}\sum_{j=0}^{k}\mu_{j}f_{j}^{p-1}\bigg)^{p^*-1}=I\!I_{i}(f)^{p^*-1}.$$ Then
$$\sup_{f\in\mathscr{F}_{I\!I}}\inf_{i\in E}I\!I_{i}(f)^{-1}\geqslant\inf_{i\in E}R_{i}(w)$$ holds for every $w\in\mathscr{W}$ and the assertion follows immediately.

\noindent{\bf Step 4}\quad Prove that $$\sup_{w\in\mathscr{W}}\inf_{i\in E}R_{i}(w)\geqslant\lambda_p.$$

If $\sum_{i\in E}\hat{\nu}_i<\infty$, then choose $f$ to be a positive function satisfying $h=fI\!I(f)^{p^*-1}<\infty.$
We have
$$h_i=\sum_{k=i}^{N}\hat{\nu}_k\bigg(\sum_{j=0}^k\mu_jf_j^{p-1}\bigg)^{p^*-1},\quad h\downdownarrows,\qquad h_i-h_{i+1}\!=\!\hat{\nu}_i\bigg(\sum_{j=0}^i\mu_jf_j^{p-1}\bigg)^{p^*-1}\!\!.$$
Let $\bar w_i=h_{i+1}/h_i$ for $i\in E$. By  a simple calculation, we obtain
$$R_i(\bar w)=\frac{-\Omega_p h(i)}{\mu_ih_i^{p-1}}=\frac{f_i^{p-1}}{h_i^{p-1}}>0.$$
If $\sum_{i\in E}\hat{\nu}_i=\infty$, then set $\bar w\equiv1$. We have $R_i(\bar w)=0$.
So
$$\sup_{w\in \mathscr{W}}\inf_{i\in E}R_i(w) \geqslant0.$$
Without loss of generality, assume that $\lambda_p>0$. By Proposition \ref{Np2}, the eigenfunction  $g$ of $\lambda_p$ is positive and strictly decreasing. Let $\bar w_i=g_{i+1}/g_i\in\mathscr{W}$. Then the assertion follows from  the fact that $R_i(\bar w)=\lambda_p$ for every $i\in E$.

\noindent{\bf Step 5} \quad We prove that the supremum in the circle arguments can be attained.

As an application of the circle arguments before Step 1, the assertion is easy in the case of $\lambda_p=0$  since
$$0=\lambda_p\geqslant\inf_{i\in E}I\!I_i(f)^{-1}\geqslant0,\qquad 0=\lambda_p\geqslant\inf_{i\in E}I_i(f)^{-1}\geqslant0$$
for every $f$ in the set defining $\lambda_p$ and
$$\lambda_p\geqslant\sup_{w\in\scr{W}}\inf_{i\in E}R_i(w)\geqslant\inf_{i\in E}R_i(\bar w)\geqslant0$$
for $\bar w$ used  in Step 4 above. In the case that $\lambda_p>0$ with eigenfunction $g$ satisfying $g_0=1$, let $\bar w_i=g_{i+1}/g_i$. Then  $R_i(\bar w)=\lambda_p$ as seen from the remarks after Proposition \ref{Np3} and
$I_i(g)=\lambda_p$ by letting $f=\lambda_p^{p^*-1}g$ in \eqref{N-D}. Moreover, we have $I\!I_i(g)=\lambda_p$ for $i\in E$ by letting $f=\lambda_p^{p^*-1}g$ in \eqref{Nf1} whenever $g_{N+1}=0$.

 It remains to rule out the probability that $g_{N+1}>0$. The Proposition \ref{Np1} below, which is proved by the variational formulas verified in Step 1, gives us the positive answer.
\begin{prop}\label{Np1} Assume that $\lambda_p>0$ and $p>1$. Let $g$ be an eigenfunction corresponding to $\lambda_p$. Then $$g_{N+1}:=\lim_{i\rightarrow N+1}g_{i}=0.$$\end{prop}
{\it Proof}\quad\rm Let $f=g-g_{N+1}$. Then $f\in\mathscr{F}_{I\!I}$. By \eqref{Nf3}, we have
$$\lambda_p^{1-p^*}f_i=\sum_{j=i}^{N}\hat{\nu}_j\bigg(\sum_{k=0}^{j}\mu_{k}{g_k}^{p-1}\bigg)^{p^*-1}.$$
Denote
$$M_i=\sum_{j=i}^{N}\hat{\nu}_j
\bigg(\sum_{k=0}^{j}\mu_{k}\bigg)^{p^*-1}.$$
If $M_i=\infty$, then
$$\lambda_p^{1-p^*}f_i=\sum_{j=i}^{N}\hat{\nu}_j\bigg(\sum_{k=0}^{j}\mu_{k}{g_k}^{p-1}\bigg)^{p^*-1}>M_ig_{N+1}.$$
There is a contradiction once $g_{N+1}\neq0$. If $M_i<\infty$, then
$$\begin{aligned}
\sup_{i\in E}I\!I_{i}(f)\!
&=\sup_{i\in E}\frac{1}{\big(g_{i}-g_{N+1}\big)^{p-1}}\bigg[\sum_{j=i}^{N}\hat{\nu}_j\bigg(\sum_{k=0}^{j}\mu_{k}{\big(g_k-g_{N+1}\big)}^{p-1}\bigg)^{p^*-1}\bigg]^{p-1}\\
&=\sup_{i\in E}\frac{1}{\lambda_p}\Bigg\{\sum_{j=i}^{N}\hat{\nu}_j
  \bigg[\sum_{k=0}^{j}\mu_{k}\big(g_k-g_{N+1}\big)^{p-1}\bigg]^{p^*-1}\\
&\qquad\qquad\quad \Bigg/\sum_{j=i}^{N}\hat{\nu}_j\bigg[\sum_{k=0}^{j}\mu_{k}g_{k}^{p-1}\bigg]^{p^*-1}\Bigg\}^{p-1}
\qquad(\text{by }  \eqref{Nf3})\\
&\leqslant\frac{1}{\lambda_p}\sup_{k\in E} \bigg(1-\frac{g_{N+1}}{g_{k}}\bigg)^{p-1}\qquad(\text{by  the proportional property})\\
&=\frac{1}{\lambda_p}\bigg(1-\frac{g_{N+1}}{g_{0}}\bigg)^{p-1}\qquad(\text{since } g\downdownarrows)
\end{aligned}$$
If $g_{N+1}>0$, then by the variational formula for lower estimates proved in Step 1 above, we have
$$\lambda_p^{-1}\leqslant\inf_{f\in\mathscr{F}_{I\!I}}\sup_{i\in E}I\!I_{i}(f)\leqslant \sup_{i\in E}I\!I_{i}(f)<\lambda_p^{-1},$$
which is a contradiction. So we must have $g_{N+1}=0$.\qquad$\Box$
\medskip

By now, we have finished the proof for the lower estimates. From this proposition, we see that the vanishing property of eigenfunction holds naturally. So the classification also holds for $N=\infty$. Combining with \eqref{Nf1}, the vanishing property also further explains where  the operator $I\!I$ comes from.
Then, we come back to the main proof of Theorem \ref{Nt1}.

For the  upper estimates,  we adopt the following circle arguments.
$$\aligned
\lambda_p
&\leqslant\inf_{f\in\widetilde{\mathscr{F}}'_{I\!I}\cup\widetilde{\mathscr{F}}_{I\!I}}\sup_{i\in {\rm supp}(f)}I\!I_{i}(f)^{-1}\\
&\leqslant \inf_{f\in\widetilde{\mathscr{F}}_{I\!I}}\sup_{i\in {\rm supp}(f)}I\!I_{i}(f)^{-1}\\&=\inf_{f\in\widetilde{\mathscr{F}}_{I}}\sup_{i\in {\rm supp}(f)}I\!I_{i}(f)^{-1}\endaligned$$
$$\aligned
&=\inf_{f\in\widetilde{\mathscr{F}}_{I}}\sup_{i\in E}{I}_{i}(f)^{-1}\\
&\leqslant\inf_{f\in\widetilde{\mathscr{F}}'_{I}}\sup_{i\in E}I_{i}(f)^{-1}\\
&\leqslant\inf_{w\in\widetilde{\mathscr{W}}}\sup_{i\in E}R_{i}(w)\\
&\leqslant\lambda_p.\endaligned$$
Since the proofs are parallel to that of the lower bounds part, we ignore most of the details here
and only mention a technique when proving
$$\inf_{f\in\widetilde{\mathscr{F}}_{I\!I}}\sup_{i\in {\rm supp}(f)}I\!I_{i}(f)^{-1}
=\inf_{f\in\widetilde{\mathscr{F}}_{I}}\sup_{i\in {\rm supp}(f)}I\!I_{i}(f)^{-1}=\inf_{f\in\widetilde{\mathscr{F}}_{I}}\sup_{i\in E}{I}_{i}(f)^{-1}.$$
To see this, we adopt a small circle arguments below:
 $$\lambda_p\leqslant\inf_{f\in\widetilde{\mathscr{F}}_{I\!I}}\sup_{i\in {\rm supp}(f)}I\!I_{i}(f)^{-1}
\leqslant\inf_{f\in\widetilde{\mathscr{F}}_{I}}\sup_{i\in {\rm supp}(f)}I\!I_{i}(f)^{-1}\leqslant\inf_{f\in\widetilde{\mathscr{F}}_{I}}\sup_{i\in E}{I}_{i}(f)^{-1}\leqslant\lambda_p.$$
The technique is about an approximating procedure, which is used to prove the last inequality above.
Recall that
$$\lambda_p^{(m)}=\inf\{D_p(f):\mu\big(|f|^p\big)=1, f=f\mathbbold{1}_{\cdot\leqslant m}\}$$
 and $\lambda_p^{(m)}\downarrow\lambda_p$ as $m\uparrow N$ (see Remark \ref{NR1}). Let $g$ be an eigenfunction of $\lambda_p^{(m)}>0$  with $g_{0}=1$.  Then   $\{g_i\}_{i=0}^{m}$ is strictly decreasing  and $g_{m+1}=0$ by letting $E$ be $E^{(m)}:=[0,m]\cap E$ in Proposition  \ref{Np3}. Extend $g$ to $E$ with $g_{i}=0$ for  $i\geqslant m+1$, we have $$g\in\widetilde{\mathscr{F}}'_{I},\qquad {\rm supp}(g)=\{0,1,\cdots,m\}\qquad \text{and }\quad \lambda_p^{(m)}=I_{k}(g)^{-1}\,\text{ for }\, k\leqslant m.$$
Hence,
$$\aligned
\lambda_p^{(m)}=\sup_{i\leqslant m}I_{i}(g)^{-1}\geqslant\inf_{f\in\widetilde{\mathscr{F}}'_{I},{\rm supp}(f)=E^{(m)}}\sup_{k\in E}I_{k}(f)^{-1}\geqslant\inf_{f\in\widetilde{\mathscr{F}}'_{I}}\sup_{k\in E}I_{k}(f)^{-1}.
\endaligned$$
Since $\widetilde{\mathscr{F}}_I'\subset\widetilde{\mathscr{F}}_I$, the right-hand side of the formula above is bounded below by $\inf_{f\in\widetilde{\mathscr{F}}_{I}}\sup_{k\in E}I_{k}(f)^{-1}$. So the required assertion follows by letting $m\rightarrow N$. $\qquad\Box$
\medskip

Instead of the approximating with finite state space used in the proof of the upper bound above, it seems more natural to use the truncating procedure for the ``eigenfunction'' $g$. However, the next result, which is easy to check by \eqref{N-D} and Proposition \ref{Np3},  shows that the procedure is not practical in general.
\begin{rem}\label{NR2}
{\rm Let $(\lambda_p,g)$ be a non-trivial  solution to eigenequation \eqref{fr2} and \eqref{B1}  with $\lambda_p>0$. Define $g^{(m)}=g\mathbbold{1}_{\leqslant m}$. Then
$$
\aligned
\min_{i\in {\rm supp}(g^{(m)})}I\!I_{i}\big(g^{(m)}\big)=\big(1-{g_{m+1}}/{g_{m}}\big)^{p-1}\big/{\lambda_p}.
\endaligned$$
In particular, the sequence $\big\{ \min_{i\in {\rm supp}(g^{(m)})}I\!I_{i}\big(g^{(m)}\big)\big\}_{m\geqslant1}$ may not converge to $\lambda_p^{-1}$ as $m\uparrow\infty.$
}\end{rem}
{\noindent \it Proof }\quad
The proof is simply an application of $f=\lambda_p^{1/(p-1)}g$ to \eqref{N-D}, based on  Proposition  \ref{Np3}.$\qquad\Box$
\medskip

 For simplicity, we write $\varphi_i=\hat{\nu}[i, N]^{p-1}$ in the proofs of  Theorems \ref{Nt2}, \ref{Nt3} and Corollary \ref{Nc1} below.

\noindent{\it Proof of Theorem $\ref{Nt2}$}\quad \rm
(a)\quad First, we prove that $\lambda_p\geqslant (k(p)\sigma_p)^{-1}$. Without loss of generality, assume that $\varphi_0<\infty$ (otherwise $ \sigma_p=\infty$). Let $f=\varphi^{1/p}<\infty$. Using the summation by parts formula, we have
$$\aligned\sum_{j=0}^{i}\mu_jf_j^{p-1}&=\mu[0,i]\varphi_i^{1/p^*}+\sum_{j=0}^{i-1}\mu[0,j]
\Big(\varphi_j^{1/p^*}-\varphi_{j+1}^{1/p^*}\Big)\\
&\leqslant\sigma_p\bigg[\varphi_i^{-1/p}+\sum_{j=0}^{i-1}\frac{1}{\varphi_j}
\Big(\varphi_j^{1/p^*}-\varphi_{j+1}^{1/p^*}\Big)\bigg]\\&\leqslant p \sigma_p\, \varphi_i^{-1/p}\endaligned$$
In the last inequality, we have used the fact that
$$\sum_{j=0}^{i-1}\frac{1}{\varphi_j}
\Big(\varphi_j^{1/p^*}-\varphi_{j+1}^{1/p^*}\Big)\leqslant (p-1)\varphi_i^{-1/p}.$$
To see this, since $\varphi_0>0$, it suffices to show that
$$\varphi_j^{1/p^*}-\varphi_{j+1}^{1/p^*}\leqslant(p-1)\varphi_j
\Big(\varphi_{j+1}^{-1/p}-\varphi_j^{-1/p}\Big).$$ Multiplying $\varphi_{j+1}^{1/p}$ on  both sides, this is  equivalent to $$p\varphi_j^{1/p^*}\varphi_{j+1}^{1/p}\leqslant(p-1)\varphi_j+\varphi_{j+1}^{1/p^*}\varphi_{j+1}^{1/p},$$
which is now obvious by Young's inequality:
$$\varphi_j^{1/p^*} \varphi_{j+1}^{1/p}\leqslant\frac{1}{p^*}
  \Big(\varphi_j^{1/p^*} \Big)^{p^*}+\frac{1}{p}\Big(\varphi_{j+1}^{1/p}\Big)^{p}.$$
Since
$$\aligned\frac{1}{\nu_i}=\hat{\nu}_i^{p-1}=\Big(\varphi_i^{p^*-1}-\varphi_{i+1}^{p^*-1}\Big)^{p-1},
\qquad \varphi_i^{{1}\big/{p(p-1)}}\varphi_{i+1}^{1/p}\leqslant \frac{1}{p}\varphi_i^{p^*-1}+\frac{1}{p^*}\varphi_{i+1}^{p^*-1},\endaligned$$
we have
\begin{eqnarray}
I_i(f)\hskip-0.5cm&&=\frac{1}{\nu_i\big(\varphi_i^{1/p}-\varphi_{i+1}^{1/p}\big)^{p-1}}\sum_{j=0}^{i}\mu_j\varphi_j^{1/p^*}
\nonumber\\
&&\leqslant
\frac{p\sigma_p\varphi_i^{-1/p}}{\nu_i\big(\varphi_i^{1/p}-\varphi_{i+1}^{1/p}\big)^{p-1}}\nonumber\\
&&=p\sigma_p\Big[\Big(\varphi_i^{p^*-1}-\varphi_{i+1}^{p^*-1}\Big)\Big/
\Big(\varphi_i^{p^*-1}-\varphi_i^{{1}/{p(p-1)}}\varphi_{i+1}^{1/p}\Big)\Big]^{p-1}\nonumber\\
&&\leqslant p{p^*}^{p-1}\sigma_p.\label{star}
\end{eqnarray}
Then the required assertion follows by Theorem \ref{Nt1}\,(1).

(b)\quad Next, we  prove that $\lambda_p\geqslant\sigma_p^{-1}$.
Let $f=\hat{\nu}[\cdot\vee n,m]\mathbbold{1}_{\cdot\leqslant m}$ for some $m, n\in E$ with $n<m$. Then $f\in\widetilde{\mathscr{F}}_{I}$ and $f_i-f_{i+1}=\hat{\nu}_i\mathbbold{1}_{n\leqslant i\leqslant m}$. By convention $1/0=\infty$, we have
$$I_i(f)=\bigg(\sum_{k=0}^n\mu_k\hat{\nu}[n,m]^{p-1}+\sum_{k=n+1}^i\mu_k\hat{\nu}[k,m]^{p-1}\bigg)\mathbbold{1}_{[n,m]}+\infty\mathbbold{1}_{[n,m]^{c}}.$$
So $$\aligned\lambda_p^{-1}&=\sup_{f\in\widetilde{\mathscr{F}}_{I}}\inf_{i\in E}I_i(f)\\
&\geqslant\inf_{i\in E}I_i(f)\\
&=\inf_{n\leqslant i\leqslant m}I_i(f)\\
&=\inf_{n\leqslant i\leqslant m}\bigg(\sum_{k=0}^n\mu_k\hat{\nu}[n,m]^{p-1}+\sum_{k=n+1}^i\mu_k\hat{\nu}[k,m]^{p-1}\bigg)\\
&=\sum_{k=0}^n\mu_k\hat{\nu}[n,m]^{p-1},\qquad m>n. \endaligned$$
The assertion that $\lambda_p^{-1}\geqslant\sigma_p$ follows by letting $m\to N$.

(c)\quad At last, if $\hat{\nu}[1,\infty)=\infty,$
then $\lambda_p=0$ is obvious.
If
$$\sum_{k=1}^{\infty}\hat{\nu}_k\mu[0,k]^{p^*-1} <\infty,$$
then $$\varphi_n\mu[0,n]
=\bigg(\sum_{k=n}^{\infty}\hat{\nu}_k\mu[0,n]^{p^*-1}\bigg)^{p-1}
\leqslant\bigg(\sum_{k=1}^{\infty}\hat{\nu}_k\mu[0,k]^{p^*-1}\bigg)^{p-1}<\infty.$$ So $\sigma_p=\infty$ and $\lambda_p=0$. $\qquad\Box$
\medskip

\noindent{\it Proof of Theorem $\ref{Nt3}$}
\quad By definitions of $\{{\bar\delta}_n\}$ and $\{\delta_n'\}$, using the proportional property, it is not hard to prove most of the results except that $\bar\delta_{n+1}\geqslant\delta_n' (n\geqslant1)$. Put $g=f_{n+1}^{(\ell, m)}$ and $f=f_n^{(\ell, m)}$. Then $g=fI\!I(f)^{p^*-1}\mathbbold{1}_{\cdot\leqslant m}$. By a simple calculation, we have
$$\big(g_{i}-g_{i+1}\big)^{p-1}=\frac{1}{\nu_i}\sum_{k=0}^{i}\mu_{k}f_{k}^{p-1},\qquad i\leqslant m.$$
Inserting this term into $D_p(g)$, we obtain
$$D_p(g)=\sum_{i=0}^m\nu_i\big(g_{i}-g_{i+1}\big)^{p-1}\big(g_{i}-g_{i+1}\big)=\sum_{i=0}^m\sum_{k=0}^{i}\mu_{k}f_{k}^{p-1}\big(g_{i}-g_{i+1}\big).$$
Noticing $g_{m+1}=0$ and exchanging  the order of the sums,  we obtain
$$D_p(g)=\sum_{k=0}^{m}\mu_{k}f_{k}^{p-1}\sum_{i=k}^{m}\big(g_{i}-g_{i+1}\big)=\sum_{k=0}^{m}\mu_{k}f_{k}^{p-1}g_{k}\leqslant \sum_{k=0}^{m}\mu_{k}g_{k}^{p}\max_{0\leqslant i\leqslant m}\frac{f_{k}^{p-1}}{g_{k}^{p-1}},$$
i.e.,
$$D_p(g)\leqslant\mu\big(|g|^{p}\big)\sup_{0\leqslant i\leqslant m}I\!I_{i}(f)^{-1}.$$
So the required assertion follows by definitions of $\bar{\delta}_{n+1}$ and $\delta_n'$.$\qquad\Box$
\medskip

Most of the results in Corollary \ref{Nc1} can be obtained from Theorem \ref{Nt3} directly. Here, we study only those assertions concerning  $\delta_1'$ and $\bar\delta_1$.

\noindent{\it Proof of Corollary $\ref{Nc1}$}\quad\rm
(a)\quad We compute  $\delta_1'$ first.
Since $p>1$ and $$\frac{1}{\hat{\nu}[i,m]}\sum_{j=i}^{m}\hat{\nu}_j\bigg(\sum_{k=0}^{j}\mu_{k}\hat{\nu}[k\vee\ell, m]^{p-1}\bigg)^{p^*-1}$$
is increasing in $i\in [\ell, m]$ (not hard to check),
we have
$$\aligned\min_{i\leqslant m}I\!I_{i}\big(f_1^{(\ell,m)}\big)
=\bigg[\frac{1}{\hat{\nu}[\ell, m]}\sum_{j=\ell}^{m}\hat{\nu}_j\bigg(\sum_{k=0}^{j}\mu_{k}
\hat{\nu}[k\vee\ell, m]^{p-1}\bigg)^{p^*-1}\bigg]^{p-1}.\endaligned$$
We claim that
$$\aligned\delta_{1}'=\sup_{\ell\in E}\frac{1}{\varphi_{\ell}}
\bigg[\sum_{j=\ell}^N\hat{\nu}_j\bigg(\sum_{k=0}^{j}\mu_{k}\varphi_{k\vee \ell}\bigg)^{p^*-1}\bigg]^{p-1}\endaligned$$
because
$$\frac{1}{\hat{\nu}[\ell, m]^{p-1}}\bigg[\sum_{j=\ell}^{m}\hat{\nu}_j\bigg(\sum_{k=0}^{j}\mu_{k}\hat{\nu}[k\vee\ell, m]^{p-1}\bigg)^{p^*-1}\bigg]^{p-1}$$
is increasing in $ m\; (m>\ell)$.
To see this, it suffices to show that
$$\aligned&\frac{1}{\hat{\nu}[\ell, m+1]}\sum_{j=\ell}^{m+1}\hat{\nu}_j\bigg(\sum_{k=0}^{j}\mu_{k}\hat{\nu}[k\vee\ell, m+1]^{p-1}\bigg)^{p^*-1}\\
&\qquad\geqslant \frac{1}{\hat{\nu}[\ell, m]}\sum_{j=\ell}^{m}\hat{\nu}_j\bigg(\sum_{k=0}^{j}\mu_{k}\hat{\nu}[k\vee\ell, m]^{p-1}\bigg)^{p^*-1}.\endaligned$$
Equivalently,
$$\sum_{j=\ell}^{m+1}\hat{\nu}_j\bigg(\sum_{k=0}^{j}\mu_{k}\frac{\hat{\nu}[k\vee\ell, m+1]^{p-1}}{\hat{\nu}[\ell, m+1]^{p-1}}\bigg)^{p^*-1}\geqslant \sum_{j=\ell}^{m}\hat{\nu}_j\bigg(\sum_{k=0}^{j}\mu_{k}\frac{\hat{\nu}[k\vee\ell, m]^{p-1}}{\hat{\nu}[\ell, m]^{p-1}}\bigg)^{p^*-1}.$$
It suffices to show that
$$\frac{\hat{\nu}[k\vee\ell, m+1]}{\hat{\nu}[\ell, m+1]}\geqslant \frac{\hat{\nu}[k\vee\ell, m]}{\hat{\nu}[\ell, m]},\qquad k\in E.$$
When $k\leqslant \ell$, the required assertion is obvious. When $k>\ell$, the inequality is just
$$\frac{\hat{\nu}[k,m+1]}{\hat{\nu}[k,m]}\geqslant \frac{\hat{\nu}[\ell, m+1]}{\hat{\nu}[\ell, m]}.$$
Noticing that $\hat{\nu}[i,m+1]=\hat{\nu}_{m+1}+\hat{\nu}[i,m]$ for any fixed $i\leqslant m$ and $\hat{\nu}[k,m]<\hat{\nu}[\ell, m]$ for $k>\ell$, we have
$$\frac{\hat{\nu}[k,m+1]}{\hat{\nu}[k,m]}=1+\frac{\hat{\nu}_{m+1}}{\hat{\nu}[k,m]}>1+\frac{\hat{\nu}_{m+1}}{\hat{\nu}[\ell, m]}
 =\frac{\hat{\nu}[\ell, m+1]}{\hat{\nu}[\ell, m]},$$
 and then the required monotone property follows.

(b)\quad Computing  $\bar{\delta}_{1}$.
Since
$$\aligned\mu\big(|f_1^{(\ell,m)}|^{p}\big)&=\sum_{j=0}^{m}\mu_{j}\hat{\nu}[\ell\vee j,m]^{p}
=\mu[0, \ell]\hat{\nu}[\ell, m]^{p}+\sum_{j=\ell+1}^{m}\mu_{j}\hat{\nu}[j,m]^{p},\endaligned$$
and
$$\aligned D_p\big(f_1^{(\ell,m)}\big)&=\sum_{j=0}^{m}\nu_j\Big(f_1^{(\ell,m)}(j)-f_1^{(\ell,m)}(j+1)\Big)^{p}\\
&=\sum_{j=\ell}^{m}\nu_j\hat{\nu}_j^p\\
&=\hat{\nu}[\ell, m]\qquad\Big(\text{since } \hat{\nu}_k=\nu_k^{1-p^*}\Big),
\endaligned$$
we have
$$\aligned
\frac{\mu\big(|f_1^{(\ell,m)}|^{p}\big)}{D_p\big(f_1^{(\ell,m)}\big)}=\hat{\nu}[\ell, m]^{p-1}\mu[0, \ell]+\frac{1}{\hat{\nu}[\ell, m]}\sum_{k=\ell+1}^m\mu_k\hat{\nu}[k,m]^p.
\endaligned$$
So
\begin{eqnarray*}
\bar{\delta}_{1}\hskip-0.7cm&&=\sup_{\ell,m\in E:\,\ell<m}\bigg(\hat{\nu}[\ell, m]^{p-1}\mu[0, \ell]+\frac{1}{\hat{\nu}[\ell, m]}\sum_{k=\ell+1}^m\mu_k\hat{\nu}[k,m]^p\bigg).
\end{eqnarray*}
The assertion on ${\bar\delta}_1$ follows immediately once we show that
$$\hat{\nu}[\ell, m]^{p-1}\mu[0, \ell]+\frac{1}{\hat{\nu}[\ell, m]}\sum_{k=\ell+1}^m\mu_k\hat{\nu}[k,m]^p$$
is increasing in $ m\; (\ell<m)$.
To see this,  it suffices to show that
$$\frac{1}{\hat{\nu}[\ell, m]}\sum_{k=\ell+1}^m\mu_k\hat{\nu}[k,m]^p\leqslant\frac{1}{\hat{\nu}[\ell, m+1]}\sum_{k=\ell+1}^{m+1}\mu_k\hat{\nu}[k,m+1]^p,$$
or equivalently,
$$\frac{\mu_{m+1}}{\hat{\nu}[\ell, m+1]}\hat{\nu}_{m+1}^{p}+\sum_{k=\ell+1}^{m}\mu_k\bigg(\frac{\hat{\nu}[k,m+1]^p }{\hat{\nu}[\ell, m+1]}-\frac{\hat{\nu}[k,m]^p}{\hat{\nu}[\ell, m]}\bigg)\geqslant0.$$
Since $p>1$ and $k>\ell$, we have $$\bigg(\frac{\hat{\nu}[k,m+1]}{\hat{\nu}[k,m]}\bigg)^p>\frac{\hat{\nu}[k,m+1]}{\hat{\nu}[k,m]}=1+\frac{\hat{\nu}_{m+1}}{\hat{\nu}[k,m]}>1+\frac{\hat{\nu}_{m+1}}{\hat{\nu}[\ell, m]}=\frac{\hat{\nu}[\ell, m+1]}{\hat{\nu}[\ell, m]}.$$
So the required assertion holds.

(c)\quad We compare $\bar\delta_1$ with $\sigma_p$ and $\delta_1'$.

 For the convenience of comparison of   ${\bar\delta}_1$ with $\sigma_p$ and $\delta_1'$,  we rewrite ${\bar\delta}_1$ as follows.
$$\aligned {\bar\delta}_1&=\sup_{l\in E}\bigg(\varphi_{\ell}\mu[0, \ell]+\frac{1}{\varphi_{\ell}^{p^*-1}}\sum_{k=\ell+1}^N\mu_k\varphi_k^{p(p^*-1)}\bigg).
\endaligned$$
By definition of $\sigma_p$, it is clear that $\bar\delta_1\geqslant\sigma_p$. To compare $\bar\delta_1$ with $\delta_1'$, we further change the form of $\bar\delta_1$. By definition of $\varphi$, we have
$$\aligned
\sum_{j=0}^{N}\mu_{j}\varphi_{j \vee m}^{p^*}&=
\sum_{j=0}^{m-1}\mu_{j}\varphi_m\sum_{k=m}^N\hat{\nu}_k+\sum_{j=m}^N\mu_{j}\varphi_j\sum_{k=j}^N\hat{\nu}_k\\
&=
\sum_{k=m}^N\hat{\nu}_k\sum_{j=0}^{m-1}\mu_{j}\varphi_m+\sum_{k=m}^N\hat{\nu}_k\sum_{j=m}^k \mu_{j}\varphi_j.\endaligned$$
So
$$\aligned
\bar{\delta}_1&=\sup_{l\in E}\bigg(\frac{1}{\varphi_{\ell}^{p^*-1}}\sum_{k=0}^N\mu_k\varphi_{k\vee \ell}^{p^*}\bigg)=\sup_{m\in E}\frac{1}{\varphi_{m}^{p^*-1}}\sum_{k=m}^N\hat{\nu}_k\sum_{j=0}^{k}\mu_{j}\varphi_{m\vee j}.\endaligned$$
Denote $a_{\ell}(k)=\hat{\nu}_k/\varphi_{\ell}^{p^*-1}$.
Then $\sum_{k=\ell}^Na_{\ell}(k)=1$
$\big($i.e., $\{a_{\ell}(k): k=\ell,\ldots, N\}$ is a probability measure on
$\{\ell, \ell+1,\ldots, N\}\big)$. By the increasing property of the moments $\mathbb{E}\big(|X|^s\big)^{1/s}$ in $s>0$, it follows that
$$\aligned
\delta_1'&=\sup_{\ell\in E}\bigg[\sum_{j=\ell}^Na_{\ell}(j)\bigg(\sum_{k=0}^{j}\mu_{k}\varphi_{k\vee \ell}\bigg)^{p^*-1}\bigg]^{p-1}\\
&\geqslant\sup_{\ell\in E}\sum_{j=\ell}^Na_{\ell}(j)\sum_{k=0}^{j}\mu_{k}\varphi_{k\vee \ell}\quad (\text{if } p^*-1>1)\\
&=\bar\delta_1. \endaligned$$
Hence, $\bar\delta_1\leqslant\delta_1'$ for $1<p\leqslant2$. Otherwise, $\bar\delta_1\leqslant\delta_1'$ for $p\geqslant2$.

(d)\quad At last, we prove that $\bar{\delta}_1\leqslant p\sigma_p$. Using the summation by parts formula,   we have
$$\aligned
\sum_{j=0}^N\mu_{j}\varphi_{j \vee \ell}^{p^*}=\sum_{j=\ell}^N
\Big(\varphi_{j}^{p^*}-\varphi_{j+1}^{p^*}\Big)\mu[0,j].
\endaligned$$
Hence,
$$\aligned\frac{1}{\varphi_{m}^{p^*-1}}\sum_{j=0}^N\mu_{j}\varphi_{j\vee m}^{p^*}
&=\frac{1}{\varphi_{m}^{p^*-1}}\sum_{j=m}^N\Big(\varphi_{j}^{p^*}-\varphi_{j+1}^{p^*}\Big)\mu [0,j]\\
&\leqslant\sigma_p\frac{1}{\varphi_{m}^{p^*-1}}\sum_{j=m}^N\frac{1}{\varphi_j}
   \Big(\varphi_{j}^{p^*}-\varphi_{j+1}^{p^*}\Big)\\
&\leqslant\sigma_p\sum_{j=m}^N\frac{1}{\varphi_j}
   \Big(\varphi_{j}^{p^*}-\varphi_{j+1}^{p^*}\Big)\bigg/\sum_{j=m}^N\Big(\varphi_{j}^{p^*-1}-\varphi_{j+1}^{p^*-1}\Big)\\
&\qquad\qquad\qquad\qquad\qquad\qquad\big(\text{since }\varphi_{N+1}=0\big).
\endaligned$$
By Young's inequality, we have $$\varphi_j\varphi_{j+1}^{p^*-1}\leqslant\frac{1}{p^*}\varphi_j^{p^*}+\frac{1}{p}\varphi_{j+1}^{p^*}.$$ Combining this inequality with the proportional property, we obtain
$$\aligned
\frac{1}{\varphi_{m}^{p^*-1}}\sum_{j=0}^N\mu_{j}\varphi_{j\vee m}^{p^*}&\leqslant\sigma_p\sup_{j\in E}\frac{\varphi_{j}^{p^*}-\varphi_{j+1}^{p^*}}{\varphi_j\big(\varphi_{j}^{p^*-1}-\varphi_{j+1}^{p^*-1}\big)}\\
&=\sigma_p\sup_{j\in E}\frac{\varphi_{j}^{p^*}-\varphi_{j+1}^{p^*}}{\varphi_{j}^{p^*}-\varphi_j\varphi_{j+1}^{p^*-1}}\\
&\leqslant p \sigma_p.
\endaligned$$
So the assertion holds .$\qquad\Box$
\noindent\\[4mm]

\setcounter{section}{4}
\setcounter{thm}{0}
\noindent{\bbb 4\quad DN-case}\\[0.1cm]
In this section, we use the same notations as the last section because they play the same role. However, they have different meaning in different sections.
Set $E=\{i\in \mathbbold{N}:1\leqslant i<N+1\}$ . Let $\{\mu_i\}_{i\in E}$ and $\{\nu_i\}_{i\in E}$ be two positive sequences.
Similar to the ND-case,  we have  the discrete version of $p$-Laplacian eigenvalue problem with DN-boundaries:
\begin{eqnarray*}
\text{`Eigenequation':}& \Omega_p\,g(k)=-\lambda\mu_k|g_k|^{p-2}g_k,\quad k\in E;\nonumber\\
\text{DN-boundary  conditions:}&g_0=0 \text{ and } g_{N+1}=g_N \text{ if } N<\infty,
\end{eqnarray*}
where $$\Omega_p\,f(k)=\nu_{k+1}|f_{k+1}-f_k|^{p-2}(f_{k+1}-f_k)-\nu_k|f_k-f_{k-1}|^{p-2}(f_k-f_{k-1}),\quad p>1,$$
$\nu_{N+1}:=0$ if $N<\infty$ and $\nu_{N+1}:=\lim_{i\to \infty}\nu_i$ if $N=\infty$.
Let $\lambda_p$ denote  the first eigenvalue. Then
\begin{equation}\label{DNV}
\lambda_p=\inf\big\{D_p(f)\big/\mu\big(|f|^p\big): f\neq0, D_p(f)<\infty\big\}.
\end{equation}
where
$$\mu\big(f\big)=\sum_{k\in E} \mu_k f_k\leqslant \infty,\quad D_{p}(f)=\sum_{k\in E}\nu_k|f_{k}-f_{k-1}|^p,\qquad f_0=0.$$
The constant $\lambda_p$ describes the optimal constant $A=\lambda_p^{-1}$ in the following {\it weighted Hardy inequality} :
$$\mu\big(|f|^p\big)\leqslant A D_p(f),\qquad f(0)=0,$$
or equivalently,
$$\|f\|_{L^p(\mu)}\leqslant A^{1/p}\|\partial^- f\|_{L^p(\nu)},\qquad f(0)=0, $$
where $\partial^- f(k)=f_{k-1}-f_{k}$.
 In other words, we are studying again the weighted Hardy inequality in this section. In view of this, by a duality \rf{KP}{Page 13}, the optimal constant $\lambda_p^{-1/p}$ in the last inequality coincides with $\lambda_{p^*}^{-1/p^*}$ which is the optimal constant in the inequality
 $$\|f\|_{L^{p^*}(\nu^{1-p^*})}\leqslant {A'}^{1/p^*}\|\partial^+ f\|_{L^{p^*}(\mu^{1-p^*})},\qquad f_{N+1}=0, $$
where $\partial^+ f(k)=f_{k+1}-f_{k}$,
 studied in Section 2. However, due to the difference of boundaries in these two cases, the variational formulas and the approximating procedure are different
 (cf. \cite{Chen1}). Therefore, it is worthy to present some details here.
  Similar notations as Section 2 are defined as follows. Define $\hat{\nu}_j=\nu_j^{1-p^*}$  for  $j\in E$, and
$$\aligned
&I_{i}(f)=\frac{1}{\nu_i(f_{i}-f_{i-1})^{p-1}}\sum_{j=i}^{N}\mu_{j}f_{j}^{p-1}\qquad(\text{single summation form}),\\
&I\!I_{i}(f)=\frac{1}{f_{i}^{p-1}}\bigg[\sum_{j=1}^{i}\hat{\nu}_j\bigg(\sum_{k=j}^{N}\mu_{k}f_{k}^{p-1}\bigg)^{p^*-1}\bigg]^{p-1}\;\quad(\text{double summation form}),\\
&R_{i}(w)\!=\!\mu_i^{-1}\big[\nu_i\big(1-w_{i-1}^{-1}\big)^{p-1}\!-\!\nu_{i+1}\big(w_{i}-1\big)^{p-1}\big],\; w_{0}\!:=\!\infty\;\;(\text{difference form}).\endaligned$$
For the lower bounds, the domains of the operators are defined respectively  as follows:
 $$\aligned
&\mathscr{F}_{I}=\{f:  f>0\; \text{and is strictly increasing on}\; E\},\\
&\mathscr{F}_{I\!I}=\{f:\ f>0\;  \text{on}\; E\},\\
&\mathscr{W}=\{w: w_{i}>1\, \text{ for } i\in E\}.
\endaligned$$
Note that the test function given in $\scr{F}_I$ is different from that given in Section 2. Again, this is due to the property of eigenfunction (which is proved in Proposition \ref{DP1} later).
For the upper bounds,  we need modify these sets as follows:
$$
\aligned
&\!\widetilde{\mathscr{F}}_{I}\!\!=\!\{f\!\!: \exists \, m\!\in\!\! E \text{ such that } f \text{ is strictly increasing on } \{1,\cdots\!,m\} \text{ and } f_{\cdot}\!\!=\!\!f_{\cdot\wedge m}\}\!,\\
&\widetilde{\mathscr{F}}_{I\!I}=\{f:f_{\cdot}=f_{\cdot\wedge m}>0\ \text{for some}\ m\in E\},\\
&\widetilde{\mathscr{W}}=\bigcup\limits_{m\in E}\Big\{w: 1< w_{i}<1+\nu_i^{p^*-1}\big(1-w_{i-1}^{-1}\big) \nu_{i+1}^{1-p^*}\text{ for }1\leqslant i\leqslant m-1 \\
&\qquad\qquad\qquad\quad\text{ and } w_{i}=1\text{ for }i\geqslant m\Big\}.
\endaligned
$$
Define $\widetilde R$ acting on $\widetilde{\mathscr{W}}$ as a modified form of $R$  by replacing $\mu_m$ with $\tilde{\mu}_m:=\sum_{k=m}^{N}\mu_k$ in $R_{i}(w)$ for the same $m$ in $\widetilde{\mathscr{W}}$.
The change of $\mu_m$ is due to the Neumann boundary at right endpoint. Note that if $w_{i}=1$ for every $i\geqslant m$, then $$\widetilde{R}_{i}(w)=R_{i}(w)=0,\qquad i>m.$$ Besides, we also need the following set:
$$\aligned
&\widetilde{\mathscr{F}}'_{I\!I}=\big\{f:f>0 \text{ on }E \text{ and }  fI\!I(f)^{p^*-1}\in L^{p}(\mu)\big\}.
\endaligned$$
If $\sum_{i\in E}{\mu_{i}}=\infty$, let $f_{i}=1$ for $i\in E$ and $f_{0}=0$. Then
 $$\aligned D_{p}(f)=\sum_{k=1}^{N}\nu_k|f_{k}-f_{k-1}|^{p}=\nu_1<\infty\quad\text{ and }\quad \mu\big(|f|^p\big)=\infty.\endaligned$$
So $\lambda_{p}=0$ by \eqref{DNV}. If $\sum_{i\in E}\mu_{i}<\infty$, then as we will prove later (Lemma \ref{D-00}) that $\lambda_p$ coincides with $$\lambda_p^{[1]}:=\inf\big\{D_{p}(f):\mu\big(|f|^{p}\big)=1\big\}.$$
Actually, the later is also coincides with
$${\lambda}_{p}^{[1]}=\inf\big\{D_{p}(f): \mu\big(|f|^{p}\big)=1, f_{i}=f_{i\wedge m}\ \mbox{for some}\ m\in E\big\}.$$

 Now we introduce the main results, many of which are parallel to that in Section 2. However, the exchange of boundary conditions `D' and `N' makes  many difference. For example, the results related to $\widetilde R$,  the definition of  $\sigma_p$ (see Theorem \ref{Dt2} later), and so on.
\begin{thm}\label{Dt1}
Assume that  $p>1$ and $\sum_{k=1}^N\mu_i<\infty$. Then the following variational formulas holds for $\lambda_p$ $(\text{equivalently, }\lambda_p^{[1]}\text{ or } \lambda_p^{[2]})$.
\begin{itemize}\setlength{\itemsep}{-0.8ex}
\item[$(1)$] Single summation forms:
$$\aligned
\inf_{f\in\widetilde{\mathscr{F}}_{I}}\sup_{i\in E}{I}_{i}(f)^{-1}=\lambda_{p}=\sup_{f\in\mathscr{F}_{I}}\inf_{i\in E}I_{i}(f)^{-1}.
\endaligned$$
\item[$(2)$] Double summation forms:
$$\aligned
&\lambda_p=\inf_{f\in\widetilde{\mathscr{F}}_{I\!I}}\sup_{i\in E}I\!I_{i}(f)^{-1}=\inf_{f\in\widetilde{\mathscr{F}}_{I}}\sup_{i\in E}I\!I_{i}(f)^{-1}=\inf_{f\in\widetilde{\mathscr{F}}'_{I\!I}\cup\widetilde{\mathscr{F}}_{I\!I}}\sup_{i\in E}I\!I_{i}(f)^{-1},\\
&\lambda_{p}=\sup_{f\in\mathscr{F}_{I\!I}}\inf_{i\in E}I\!I_{i}(f)^{-1}=\sup_{f\in\mathscr{F}_{I}}\inf_{i\in E}I\!I_{i}(f)^{-1}.
\endaligned$$
\item[$(3)$]Difference forms:
$$\aligned
\inf_{w\in\widetilde{\mathscr{W}}}\sup_{i\in E}\widetilde R_{i}(w)=\lambda_{p}=\sup_{w\in\mathscr{W}}\inf_{i\in E}R_{i}(w).
\endaligned$$
\end{itemize}
\end{thm}

 As an application of the variational formulas in Theorem \ref{Dt1}\,(1), we have the following theorem. This result was   known in 1990's
 (cf. \rf{KMP}{Page 58, Theorem7} plus the duality technique, cf. \rf{KP}{Page 13}). See also \cite{Mao}. It can be regarded as a dual of Theorem \ref{Nt2}.
\begin{thm}\label{Dt2}
$(\text{\rm Basic estimates})$ For $p>1$, we have  $\lambda_{p}$ $($or equivalently, $\lambda_p^{[1]}$ or $\lambda_p^{[2]}$ provided $\sum_{k\in E}\mu_k<\infty$$)$ is positive  if and only if
$\sigma_p<\infty$, where
 $$\sigma_p=\sup_{n\in E}\big(\mu[n,N]\hat{\nu}[1,n]^{p-1}\big).$$
More precisely,
$$\big(k(p)\sigma_p\big)^{-1}\leqslant\lambda_{p}\leqslant\sigma_p^{-1},$$
where $k(p)=p{p^*}^{p-1}$.
In particular, we have $\lambda_p=0$
if $\sum_{i\in E}\mu_{i}=\infty$ and $\lambda_{p}>0$ if $N<\infty$, or $\sum_{k=1}^{\infty}\mu[k, N]^{p^*-1}\nu_k <\infty$, or $\sum_{k=1}^{\infty}\big(\mu_k+\hat{\nu}_k\big)<\infty$.
\end{thm}

The next result is an application of the variational formulas in Theorem \ref{Dt1}\,(2). It is interesting that the result is not a direct dual of Theorem \ref{Nt3}.
\begin{thm}\label{Dt3}
 $(\text{\rm Approximating  procedure})$ Assume that $p>1$, $\sum_{k\in E}\mu_k<\infty$ and $\sigma_p<\infty$. Then the following assertions hold.
\begin{itemize}\setlength{\itemsep}{-0.8ex}
\item[$(1)$] Define
$$f_{1}=\hat{\nu}[1, \cdot]^{1/p^*},\;\; f_{n}=f_{n-1}I\!I(f_{n-1})^{p^*-1}\,(n\geqslant 2),\quad \delta_{n}=\sup_{i\in E}I\!I_{i}(f_{n}).$$
Then $\delta_{n}$ is decreasing in $n$ and $$\aligned
\lambda_{p}\geqslant\delta_{\infty}^{-1}\geqslant\cdots\geqslant\delta_{1}^{-1}\geqslant(k(p)\sigma_p)^{-1}.
\endaligned$$
\item[$(2)$] For fixed $m\in E$, define
$$\aligned
f_{1}^{(m)}=\hat{\nu}[1, \cdot\wedge m],\qquad f_{n}^{(m)}=f_{n-1}^{(m)}I\!I\big(f_{n-1}^{(m)}\big)(\cdot\wedge m)^{p^*-1}, \quad n\geqslant2,
\endaligned$$ and
$$\delta_n'=\sup_{m\in E}\inf_{i\in E}I\!I_{i}(f_{n}^{(m)}).$$
Then $\delta_{n}'$ is increasing in $n$ and
$$\aligned
\sigma_p^{-1}\geqslant{\delta_{1}'}^{-1}\geqslant\cdots\geqslant{\delta_{\infty}'}^{-1}\geqslant\lambda_{p}.
\endaligned$$
Next, define
$$\aligned\bar{\delta}_n=\sup_{m\in E}\frac{\mu\big({f_{n}^{(m)}}^p\big)}{D_p\big(f_{n}^{(m)}\big)},\qquad
 n\in E.\endaligned$$
 Then $\bar{\delta}_{n}\geqslant\lambda_{p}$ and $\bar{\delta}_{n+1}\geqslant\delta_{n}'$ for every $n\geqslant1$.
\end{itemize}
 \end{thm}
\begin{cor}$(\text{\rm Improved estimates})$\label{Dcor1}
Assume that  $\sum_{k\in E}\mu_k<\infty$. For $p>1$, we have
$$\sigma_p^{-1}\geqslant{\delta_1'}^{-1}\geqslant\lambda_p^{-1}\geqslant\delta_1^{-1}\geqslant(k(p)\sigma_p)^{-1},$$
where
$$\aligned
\delta_1&=\sup_{i\in E}\bigg[\hat{\nu}[1, i]^{-1/p^*}\sum_{j=1}^i\hat{\nu}_j\bigg(\sum_{k=j}^N\mu_k\hat{\nu}[1, k]^{(p-1)/p^*}\bigg)^{p^*-1}\bigg]^{p-1},\\
\delta_1'&=\sup_{m\in E}\frac{1}{\hat{\nu}[1, m]^{p-1}}\bigg[\sum_{j=1}^m\hat{\nu}_j\bigg(\sum_{k=j}^{N}\mu_{k}\hat{\nu}[1, k\wedge m]^{p-1}\bigg)^{p^*-1}\bigg]^{p-1}.\endaligned$$
Moreover,
$$\bar\delta_{1}=\sup_{m\in E}\frac{1}{\hat{\nu}[1, m]}\sum_{j=1}^{N}\mu_{j}\,\hat{\nu}[1, j\wedge m]^p\in[\sigma_p, p\sigma_p],$$
and $\bar\delta_1\geqslant\delta_1'$ for $p\geqslant2$, $\bar\delta_1\leqslant\delta_1'$ for $1<p\leqslant2$.
\end{cor}

When $p=2$, the result that $\delta_1'=\bar\delta_1$ is also known (see \rf{Chen1}{Theorem 4.3}).
\noindent\\[4mm]

\noindent{\bbb 4.1\quad Partial proofs of main results}\\[0.1cm]
Before moving to the  proofs of the main results, we give some more descriptions of $\lambda_p$.  Define
$$\lambda_{p}^{(m)}=\inf\big\{ D_p(f): \mu\big(|f|^{p}\big)=1, f_{i}=f_{i\wedge m},\; i\in E\big\}.$$
Let
 $$\widetilde D_p(f)=\sum_{i=1}^{m}\tilde{\nu}_i|f_{i}-f_{i-1}|^p,
 \quad \quad\tilde\mu(f)=\sum_{i=1}^{m}\tilde{\mu}_{i}|f_{i}|^{p},$$
where $\tilde\nu$ and $\tilde\mu$ are defined  as follows:
$$\aligned\tilde{\nu}_i=\nu_i \quad\text{ for } i\leqslant m;\qquad  \tilde{\mu}_i=\mu_i\quad \text{ for }\; i\leqslant m-1, \qquad \tilde{\mu}_{m}=\sum_{k=m}^{N}\mu_{k}.
\endaligned$$
For  $f=f_{\cdot\wedge m}$, we have
 $${\widetilde D}_p(f)= D_p(f),\qquad \tilde\mu\big(|f|^{p}\big)=\mu\big(|f|^{p}\big).$$
 So $\lambda_p^{(m)}$ is the first eigenvalue of the local Dirichlet form $\big(\widetilde D, \scr{D}(\widetilde D)\big)$ on $E^{(m)}:=\{1,2,\cdots,m\}$ with  reflecting (Neumann) boundary at $m+1$  and absorbing (Dirichlet) boundary at $0$. Furthermore, we have the following fact.
\begin{lem}\label{D-00}
If $\sum_{i\in E}\mu_{i}<\infty$, then $\lambda_p=\lambda_p^{[1]}={\lambda}_{p}^{[2]}$. Moreover,
  $\lambda_{p}^{(m)}\downarrow\lambda_p^{[2]}$ as $m\to N$.
\end{lem}
{\it Proof}\quad Since each $f$ with $\mu(|f|^p)=\infty$ can be approximated by $f_i^{(m)}=f_{i\wedge m}$ ($m\in E$) with respect to norm $\|\cdot\|^p=D_p(\cdot)+\mu\big(|\cdot|^p\big)$, it is clear that $\lambda_p=\lambda_p^{[1]}$.

 We now  prove that $\lambda_p^{[1]}=\lambda_p^{[2]}$. It is  clear that $\lambda_p^{[2]}\geqslant\lambda_p^{[1]}$ since $\sum_{k\in E}\mu_k<\infty$. For any fixed $\varepsilon>0$, there exists $f\in L^p(\mu)$ such that $$D_p(f)/\mu\big(|f|^p\big)\leqslant \lambda_p^{[1]}+\varepsilon.$$ Let $f^{(n)}=f_{\cdot\wedge n}$. Then
  $$D_p(f^{(n)})\to D_p(f),\quad \mu\big(|f^{(n)}|^p\big)\to \mu\big(|f|^p\big)\qquad \text{ as } n\to N.$$
By definitions of $\lambda_p^{(n)}$ and $\lambda_p^{[2]}$, for large enough $n\in E$, we have
$$\lambda_p^{[2]}\leqslant\lambda_p^{(n)}\leqslant\frac{D_p(f^{(n)})}{\mu\big(|f^{(n)}|^p\big)}\leqslant\frac{D_p(f)}{\mu\big(|f|^p\big)}+\varepsilon
\leqslant\lambda_p^{[1]}+2\varepsilon\leqslant\lambda_p^{[2]}+2\varepsilon.$$
 Hence, $\lambda_p^{[1]}=\lambda_p^{[2]}$ and
  $\lambda_{p}^{(m)}\downarrow\lambda_p^{[2]}$.
\medskip

\noindent{\it Proof of Theorem $\ref{Dt1}$}\quad In parallel to the ND-case, we also adopt two circle arguments to  prove the theorem. For instance, the circle argument below is adopted  for the upper estimates:
$$\aligned
\lambda_{p}&\leqslant\inf_{f\in\widetilde{\mathscr{F}}'_{I\!I}\cup\widetilde{\mathscr{F}}_{I\!I}}\sup_{i\in E}I\!I_{i}(f)^{-1}\\
&\leqslant \inf_{f\in\widetilde{\mathscr{F}}_{I\!I}}\sup_{i\in E}I\!I_{i}(f)^{-1}\\&=\inf_{f\in\widetilde{\mathscr{F}}_{I}}\sup_{i\in E}I\!I_{i}(f)^{-1}\\
&=\inf_{f\in\widetilde{\mathscr{F}}_{I}}\sup_{i\in E}{I}_{i}(f)^{-1}\\
&\leqslant\inf_{w\in\widetilde{\mathscr{W}}}\sup_{i\in E}\widetilde R_{i}(w)\\
&\leqslant\lambda_{p}
.\endaligned$$
The proofs are similar to  the  ND-case. Here, we  present the  proofs  of the assertions  related to the operator $\widetilde R$ in the above circle argument, which are  obvious different from that in Section 2 due to the boundary conditions.

(1)\quad We  first prove that
$$\inf_{f\in\widetilde{\mathscr{F}}_{I\!I}}\sup\limits_{i\in E}I\!I_{i}(f)^{-1}
\leqslant\inf_{w\in\widetilde{\mathscr{W}}}\sup\limits_{i\in E}\widetilde R_{i}(w).$$
For
$w\in\widetilde{\mathscr{W}}$,  it is easy to check that $\widetilde R_{i}(w)>0$. Let
$u_{0}=0$ and  $u_{i}=u_{i\wedge m}>0$ for $i\in E$ such that $w_{i}(=u_{i+1}/u_i)\in \widetilde{\scr{W}}$.
Then $u_{i}$ is strictly increasing on $[0,m]$.
Put
\begin{eqnarray*}
f_{i}^{p-1}=
\left\{
\begin{array}{ll}\mu_i^{-1}
\big[\nu_i(u_i-u_{i-1})^{p-1}-\nu_{i+1}(u_{i+1}-u_{i})^{p-1}\big],&i\leqslant m-1, \\ \tilde\mu_m^{-1}\nu_m(u_{m}-u_{m-1})^{p-1},&i\geqslant m.
\end{array}
\right.
\end{eqnarray*}
We have $f\in\widetilde{\mathscr{F}}_{I\!I}$ and
 $$\aligned\Omega_p u(k)=-\mu_k f_k^{p-1},\qquad k\leqslant m-1,\endaligned$$
$$\aligned\nu_m(u_{m}-u_{m-1})^{p-1}=\sum_{j=m}^{N}\mu_{j}f_{m}^{p-1}=\sum_{j=m}^{N}\mu_{j}f_{j}^{p-1}.\endaligned$$
By a simple reorganization  and  making  summation from $1$ to $i\,(\leqslant m)$ with respect to $k$, we obtain
 $$\aligned
\sum_{k=1}^{i}\hat{\nu}_k\bigg(\sum_{j=k}^{N}\mu_{j}f_{j}^{p-1}\bigg)^{p^*-1}=u_{i}-u_{0}=u_{i},
\quad i\leqslant m.\endaligned$$
Therefore,
$$
\aligned
\widetilde R_{i}(w)^{p^*-1}=\frac{f_{i}}{u_{i}}=f_{i}\bigg[\sum_{k=1}^{i}\hat{\nu}_k\bigg(\sum_{j=k}^N\mu_{j}f_{j}^{p-1}\bigg)^{p^*-1}\bigg]^{-1}=I\!I_{i}(f)^{1-p^*},\qquad i\leqslant    m,
\endaligned$$
and then
$$
\aligned
\sup_{i\in E}\widetilde R_{i}(w)\geqslant\max_{i\leqslant  m}\widetilde R_{i}(w)\geqslant\sup_{i\in E}I\!I_{i}(f)^{-1}\geqslant\inf_{f\in\widetilde{\mathscr{F}}_{I\!I}}\sup_{i\in {\supp}(f)}I\!I_{i}(f)^{-1}.
\endaligned$$

(2)\quad We prove that $$\inf\limits_{w\in\widetilde{\mathscr{W}}}\sup\limits_{i\in E}\widetilde R_{i}(w)\leqslant\lambda_{p}.$$
Let $g$ denote the solution to the equation
$$-\Omega_pg(i)=\lambda_{p}^{(m)}\tilde\mu_i|g_i|^{p-2}g_i,\;\; g_{0}=0,\;\;
g_{m+1}:=g_m, \;\; i\in E^{(m)}:=\{0,1,\ldots,m\}.$$
Without loss of generality, assume  $g_1>0$. Then $g$ is strictly increasing (by Proposition \ref{DP1} below) and
$$\aligned
&-\nu_{i+1}\big(g_{i+1}-g_i\big)^{p-1}+\nu_i\big(g_i-g_{i-1}\big)^{p-1}=\lambda_{p}^{(m)}\mu_ig_{i}^{p-1}, \qquad i\leqslant    m-1.
\endaligned$$
$$\aligned
\nu_m\big(g_m-g_{m-1}\big)^{p-1}=\lambda_{p}^{(m)}\tilde\mu_mg_{m}^{p-1}.
\endaligned$$
That is,
$$\aligned&\nu_i\bigg(1-\frac{g_{i-1}}{g_{i}}\bigg)^{p-1}-\nu_{i+1}\bigg(\frac{g_{i+1}}{g_{i}}-1\bigg)^{p-1}=\lambda_{p}^{(m)}\mu_i,\qquad i\leqslant m-1;
\endaligned$$
$$\aligned\nu_m\bigg(1-\frac{g_{m-1}}{g_{m}}\bigg)^{p-1}=\lambda_{p}^{(m)}\tilde\mu_m.\endaligned$$
Let $\bar{w}_{i}={g_{i+1}}/{g_{i}}$ for $i\leqslant m-1$ and $\bar{w}_{i}=1$ for $i\geqslant m.$ Then $\bar{w}\in\widetilde{\mathscr{W}}$  and $\widetilde R_{i}(\bar w)=\lambda_{p}^{(m)}$ for every $i\leqslant m$.
 Therefore,
$$\aligned
\lambda_{p}^{(m)}&=\max_{0\leqslant i\leqslant m}\widetilde R_{i}(\bar w)\\
&\geqslant\inf_{w\in\widetilde{\mathscr{W}}:  w_{i}=w_{i\wedge m}}\sup\limits_{0\leqslant i\leqslant m}\widetilde  R_{i}(w)\\
&\geqslant\inf_{w\in\widetilde{\mathscr{W}}:  w_{i}=w_{i\wedge n}\text{ for some }  n\in E}\sup\limits_{i\in E}\widetilde R_{i}(w)\\
&\geqslant\inf_{w\in\widetilde{\mathscr{W}}}\sup\limits_{i\in E}\widetilde  R_{i}(w).
\endaligned$$
We obtain the required assertion by letting $m\rightarrow N$.
$\qquad\Box$

Noticing the difference between the ND- and the DN-cases, one may finish the proofs of other theorems in this section without much difficulties following Section 2 or \rf{Chen1}{Section 4}. So we ignore the details here but present one proposition below, which is essential to our study and is used to verify the last inequalities related to $R$ or $\widetilde R$ in the two circle arguments.  The proposition, whose proof is independent of the other assertions in this paper, provides the basis for imitating  the eigenfunction to construct the corresponding test functions of the operators.
\begin{prop}\label{DP1}
Assume that $g$ is a nontrivial solution to $p$-Laplacian problem with DN-boundary conditions.  Then $g$ is monotone. Moreover, $g$ is increasing provided $g_{1}>0.$\end{prop}
\rm
{\noindent \it Proof }\quad The proof is parallel to that in \rf{Chen2}{Proposition 3.4}. We give the skeleton of the proof. The proof is quite easy in the case of  $\lambda_p=0$.
For the case that $\lambda_p>0$, suppose that there exists $n\in E$ such that $g_{0}<g_{1}<\cdots<g_{n}\geqslant  g_{n+1}$. Then define $\bar{g}_{i}=g_{i\wedge n}$. By a simple calculation and \eqref{DNV}, we obtain
$$\aligned
\lambda_{p}\leqslant\frac{D_{p}(\bar{g})}{\mu\big(|\bar{g}|^{p}\big)}
=\frac{\lambda_{p}\sum_{k=1}^{n-1}\mu_{k}|g_{k}|^{p}+\nu_n|g_{n}-g_{n-1}|^{p-2}\big(g_{n}-g_{n-1}\big)g_{n}}
{\sum_{k=1}^{n-1}\mu_{k}|g_{k}|^{p}+|g_{n}|^{p}\sum_{k=n}^{N}\mu_{k}}<\lambda_{p}.\endaligned$$
In the last inequality, we have used the following fact:
$$\aligned
\nu_n|g_{n}\!-\!g_{n-1}|^{p-2}\big(g_{n}\!-\!g_{n-1}\big)g_{n}\leqslant -g(n)\Omega_pg(n)=\lambda_{p}\mu_{n}|g_{n}|^{p}<\lambda_{p}|g_{n}|^{p}\sum_{k=n}^{N}\mu_{k}
\endaligned$$
for $n<N$. Therefore, there is a contradiction and so the required assertion holds.
$\qquad\Box$
\noindent\\[4mm]

\noindent\bf{\footnotesize Acknowledgements}\quad\rm
{\footnotesize The authors thank Professors Yonghua Mao and Yutao Ma for their helpful comments and suggestions.
The work was supported in part by NSFC (Grant No.11131003), SRFDP (Grant No. 20100003110005), the ``985'' project from the Ministry of Education in China, and the Fundamental Research Funds for the Central Universities.
}\\[4mm]

\noindent{\bbb{References}}

\vspace{2mm}

\noindent {\small School of Mathematical Sciences, Beijing Normal University;}
\\ {\small Laboratory of Mathematics and Complex Systems (Beijing Normal University),}
\\ {\text{\hskip2em}\small Ministry of Education;}
\\ {\small  Beijing 100875, China}
\\{\small Emails: M.F. Chen\quad mfchen@bnu.edu.cn
\qquad Y.H. Zhang\quad zhangyh@bnu.edu.cn}
\end{document}